\documentclass{siamltex}
\usepackage{yuyuan}
\usepackage[utf8]{inputenc}
\usepackage{algorithm, algorithmic}
\usepackage[colorlinks=true, citecolor=blue]{hyperref}
\usepackage{fullpage}

\def\cHrt{\ensuremath{{\cal H}(r_t, \zeta_{2t-1})}}
\def\cHwt{\ensuremath{{\cal H}(w\tn, \zeta_{2t})}}
\def\cGt{\ensuremath{{\cal G}(w\md, \xi_t)}}
\def\aGi{\frac{\alpha_i}{\Gamma_i}}

\title{Accelerated Schemes For A Class of Variational Inequalities}
\author{Yunmei Chen\thanks{Department of Mathematics, University of Florida ({\tt yun@math.ufl.edu}). This author was partially supported by
    NSF grants DMS-1115568, IIP-1237814 and DMS-1319050.} \and
Guanghui Lan\thanks{Department of Industrial and System Engineering, University of Florida ({\tt glan@ise.ufl.edu}).
This author was partially supported by NSF grant CMMI-1000347, ONR grant N00014-13-1-0036, NSF DMS-1319050, and NSF CAREER Award CMMI-1254446.} \and
Yuyuan Ouyang\thanks{Department of Industrial and System Engineering, University of Florida ({\tt ouyang@ufl.edu})}. Part of the research was done while the author was a PhD student at the Department of Mathematics, University of Florida. This author was partially supported by AFRL Mathematical Modeling Optimization Institute.}

\begin{document}

\maketitle

\begin{abstract}
We propose a novel method, namely the accelerated mirror-prox (AMP) method, for computing the weak solutions of a class of deterministic and stochastic monotone variational inequalities (VI). The main idea of this algorithm is to incorporate a multi-step acceleration scheme into the mirror-prox method. For both deterministic and stochastic VIs, the developed AMP method computes the weak solutions with optimal rate of convergence. In particular, if the monotone operator in VI consists of the gradient of a smooth function, the rate of convergence of the AMP method can be accelerated in terms of its dependence on the Lipschitz constant of the smooth function. For VIs with bounded feasible sets, the estimate of the rate of convergence of the AMP method depends on the diameter of the feasible set. For unbounded VIs, we adopt the modified gap function introduced by Monteiro and Svaiter for solving monotone inclusion, and demonstrate that the rate of convergence of the AMP method depends on the distance from the initial point to the set of strong solutions.

\end{abstract}

\section{Introduction}
Let ${\cal E}$ be a finite dimensional vector space with inner product $\langle\cdot,\cdot\rangle$ and norm $\|\cdot\|$, and $Z$ be a non-empty closed convex set in $\cE$. Our problem of interest is to find $u^*\in Z$ that solves the following variational inequality (VI) problem:
\begin{align}
	\label{eqnProblem}
	& \langle F(u), u^* - u\rangle \leq 0, \forall u\in Z,
\end{align}
where $F$ is defined by
\begin{align}
	\label{eqnF}
	F(u) = \nabla G(u) + H(u) + J'(u),\ \forall u\in Z.
\end{align} 
In \eqref{eqnF}, $G(\cdot)$ is a general continuously differentiable function whose gradient is Lipschitz continuous with constant $L_G$, i.e.,
	\begin{equation}\label{eqnGAssumption}
		0\leq G(w)-G(v)-\langle\nabla G(w), w-v\rangle \leq \frac {L_G}2 \|w-v\|^2, \forall w, v\in Z,
	\end{equation}
$H:Z\to \cE$ is a monotone operator with Lipschitz constant $L_H$, that is, for all $w, v\in Z$,
\begin{align}
	& \langle H(w) - H(v), w - v\rangle \geq 0,\text{ and }
	\|H(w)-H(v)\|_*\leq L_H\|w-v\|,
\end{align}	
and $J'(u)\in\partial J(u)$, where $J(\cdot)$ is a relatively simple and convex function.
We denote problem \eqref{eqnProblem} by $VI(Z;G,H,J)$ or simply $VI(Z;F)$.

Observe that $u^*$ given by \eqref{eqnProblem} is often called a weak solution of $VI(Z;F)$. A related notion is a strong solution of VI. More specifically, we say that $u^*$ is a strong solution of $VI(Z;F)$ if it satisfies
\begin{align}
	\label{eqnSVI}
	\ & \langle F(u^*), u^* - u\rangle \leq 0, \forall u\in Z.
\end{align}
For any monotone operator $F$, it is well-known that strong solutions of $VI(Z,F)$ are also weak solutions, and the reverse is also true under mild assumptions (e.g., when $F$ is continuous). For example, for $F$ in \eqref{eqnF}, if $J=0$, then the weak and strong solutions of $VI(Z; G, H, 0)$ are equivalent. 

The main goal of this paper is to develop efficient solution methods for solving two types of VIs, i.e., deterministic VIs with exact information about the operator $F$, and stochastic VIs where the operator $F$ contains some stochastic components (e.g., $\nabla G$ and $H$) that cannot be evaluated exactly. We start by reviewing some existing methods for solving both these types of problems. 

\subsection{Deterministic VI}
\label{secIntroVID}

VI provides a unified framework for optimization, equilibrium and complementarity problems, and thus has been the focus of many algorithmic studies (see, e.g, \cite{korpelevich1976extragradient,rockafellar1976monotone,chen1999homotopy,nesterov1999homogeneous,solodov1999hybrid,solodov2000inexact,nemirovski2005prox,nesterov2007dual,monteiro2010complexity,juditsky2011solving}). In particular, classical algorithms for VI include, but not limited to, the gradient projection method (e.g., \cite{sibony1970methodes,bertsekas1999nonlinear}), Korpelevich's extragradient method \cite{korpelevich1976extragradient}, and the proximal point algorithm (e.g., \cite{martinet1970regularisation,rockafellar1976monotone}), etc. (see \cite{facchinei2003finite} for an extensive review and bibliography). While these
earlier studies on VI solution methods focused on their asymptotic convergence behavior (see, e.g., \cite{solodov1999new,sun1995new,tseng2000modified}), 
much recent research effort has been devoted to algorithms exhibiting strong performance guarantees in a finite number of iterations (a.k.a., iteration complexity) 
 \cite{nesterov1999homogeneous,ben-tal2005non,nesterov2007dual,nesterov2009primal,nemirovski2009robust,monteiro2010complexity,dang2012convergence}.
More specifically, Nemirovski in a seminal work~\cite{nemirovski2005prox} presented a mirror-prox method
by properly modifying Korpelevich's algorithm~\cite{korpelevich1983extrapolation}  
and show that it can achieve an ${\cal O}(1/\epsilon)$ complexity bound 
for solving VI problems with Lipschitz continuous operators (i.e., smooth VI denoted by
$VI(Z;0,H,0)$). Here $\epsilon > 0$ denotes the target accuracy in terms of a weak solution. This bound significantly improves the ${\cal O} (1/\epsilon^2)$ bound
for solving VI problems with bounded operators (i.e., nonsmooth VI) (e.g., \cite{ben-tal2005non}). Nemirovski's algorithm was further 
generalized by Auslender and Teboulle \cite{auslender2005interior} through the incorporation of a wider class of distance generating functions.
Nesterov~\cite{nesterov2007dual} has also developed a dual extrapolation method for solving smooth VI which possesses the same complexity bound
as in \cite{nemirovski2005prox}.
More recently, Monteiro and Svaiter~\cite{monteiro2010complexity} showed that the
the hybrid proximal extragradient (HPE) method \cite{solodov1999hybrida}, which covers Korpelevich's algorithm
as a special case, can also achieve the aforementioned ${\cal O}(1/\epsilon)$ complexity. Moreover,
they developed novel termination criterion for VI problems with possibly unbounded feasible set $Z$,
and derived the iteration complexity associated with HPE for solving unbounded VI problems accordingly.
Monteiro and Svaiter~\cite{monteiro2011complexity}
have also generalized the aforementioned ${\cal O} (1/\epsilon)$ complexity result for solving VI problems containing a simple nonsmooth
component (i.e., $VI(Z;0,H,J)$).

It should be noted, however, that the aforementioned studies in existing literature do not explore the fact that
the operator $F$ consists of a gradient component $\nabla G$ (see \eqref{eqnF}). 
As a result, the iteration complexity associated with any of these algorithms, when applied to a smooth convex optimization 
problem (i.e., $VI(Z;G,0,0)$), is given by ${\cal O} (1/\epsilon)$, which is
significantly worse than the well-known ${\cal O}(1/\sqrt{\epsilon})$ optimal complexity for smooth optimization \cite{nesterov1983method}. An important motivating question for our study is
whether one can utilize such structural properties of $F$ in order to further improve the efficiency of
VI solution methods. More specifically, we can easily see that the total number of gradient and/or operator evaluations 
for solving $VI(Z;G,H,J)$ cannot be smaller than
\begin{equation}
	\label{eqnOptRate}
	\cO\left(\sqrt{\frac{L_G}{\epsilon}}+\frac{L_H}{\epsilon} \right).
\end{equation}
Such a lower complexity bound is derived based on the following two observations:
\begin{enumerate}
	\renewcommand{\theenumi}{\alph{enumi})}
	\item If $H=0$, $VI(Z;G,0,0)$ is equivalent to a smooth optimization problem $\min_{u\in Z}G(u)$, and the complexity for minimizing $G(u)$ cannot be better than $\cO(\sqrt{L_G/\epsilon})$ \cite{nesterov1983method, nesterov2004introductory};
	\item If $G=0$, the complexity for solving $VI(Z; 0, H, 0)$ cannot be better than $\cO(L_H/\epsilon)$ \cite{nemirovski1992information} (see also the discussions in Section 5 of \cite{nemirovski2005prox}).
\end{enumerate}
However, 
the best-known so-far iteration complexity bound for solving $VI(Z;G,H,J)$ is given by \cite{juditsky2011solving,monteiro2010complexity},
where one needs to run these algorithms
\begin{equation}
	\label{eqnRateVI}
	\cO\left(\frac{L_G+L_H}{\epsilon} \right),
\end{equation}
iterations to  compute a weak solution of $VI(Z;G,H,J)$, and each iteration requires the computation of
both $\nabla G$ and $H$. 
It is worth noting that better iteration complexity bound has been achieved for a certain special case
of $VI(Z;G,H,J)$ where the operator $H$ is linear. In this case, Nesterov~\cite{nesterov2005smooth} showed that, by using a novel smoothing
technique, the total number of first-order iterations (i.e., iterations requiring the computation of $\nabla G$, the linear operators $H$ and its conjugate $H^*$) 
for solving $VI(Z;G,H,J)$ can be bounded by \eqref{eqnOptRate}. 
This bound has also been obtained by 
applying an accelerated primal-dual method recently developed by Chen, Lan and Ouyang~\cite{chen2013optimal}. 
Observe that the bound in \eqref{eqnOptRate} is significantly better than the one in \eqref{eqnRateVI} in terms of its dependence on $L_G$.
However, it is unclear whether similar iteration complexity bounds to those in \cite{nesterov2005smooth,chen2013optimal} can be achieved for the more general case 
when $H$ is Lipschitz continuous.

\subsection{Stochastic VI}
While deterministic VIs had been intensively investigated in the literature, the study of stochastic VIs is still quite limited. In the stochastic setting, we assume that there exists {\it stochastic oracles} ${\cal SO}_G$ and ${\cal SO}_H$ that provide unbiased estimates to the operators $\nabla G(u)$ and $H(u)$ for any test point $u\in Z$. More specifically, we assume that at the $i$-th call of ${\cal SO}_G$ and ${\cal SO}_H$ with input $z\in Z$, the oracles ${\cal SO}_G$ and ${\cal SO}_H$ output stochastic first-order information $\cG(z, \xi_i)$ and $\cH(z, \zeta_i)$ respectively, such that $\E[{\cal G}{(x, \xi_i)}] = \nabla G(x), \E[{\cal H}{(x, \zeta_i)}] = H(x)$, and
\begin{enumerate}
	\newcounter{enuAssumptions}
	\setcounter{enumi}{\theenuAssumptions}
	\renewcommand{\theenumi}{\textbf{A\arabic{enumi}}}
	\item
	\label{itmVB}
	 $\E\left[\left\|{\cal G}{(x, \xi_i)} - \nabla G(x)\right\|_*^2\right] \leq \sigma_G^2,\  \E\left[\left\|{\cal H}{(x, \zeta_i)} - H(x)\right\|_*^2\right] \leq \sigma_H^2,$
	\setcounter{enuAssumptions}{\value{enumi}}	
\end{enumerate}
where $\xi_i\in \Xi$, $\zeta_i\in \Xi$ are independently distributed random variables. It should be noted that deterministic VIs are special cases of stochastic VIs with $\sigma_G=\sigma_H=0$. To distinguish stochastic VIs from their deterministic counterparts, we will use $SVI(Z;G, H, J)$ or simply $SVI(Z;F)$ to denote problem \eqref{eqnProblem} under the aforementioned stochastic settings. 

Following the discussion around \eqref{eqnOptRate} and the complexity theory for stochastic optimization \cite{nemirovski1983problem, juditsky2011solving}, the total number of  gradient and operator evaluations for solving stochastic VI cannot be smaller than
\begin{equation}
	\label{eqnOptRateS}
	\cO\left(\sqrt{\frac{L_G}{\epsilon}} + \frac{L_H}{\epsilon} + \frac{\sigma_G^2+\sigma_H^2}{\epsilon^2}\right).
\end{equation}
The best known complexity bound for computing $SVI(Z;G,H,0)$ is given by the stochastic mirror-prox method in \cite{juditsky2011solving}. This method requires
\begin{equation} \label{eqnBestRate}
\cO\left(\frac{L_G+L_H}{\epsilon} + \frac{\sigma_G^2+\sigma_H^2}{\epsilon^2}\right)
\end{equation}
iterations to achieve the target accuracy $\epsilon > 0$ in terms of a weak solution, and each iteration requires the calls to ${\cal SO}_G$ and ${\cal SO}_H$.
Similar to the deterministic case, the above complexity bound has been improved
for some special cases, e.g., when $H=0$ or $H$ is linear. 
In particular, when $H=0$, $SVI(Z,F)$ is equivalent to the stochastic minimization problem  of$\min_{u\in Z} G(u) + J(u)$, 
Lan first presented in \cite{lan2012optimal} (see more general results in \cite{ghadimi2012optimal,ghadimi2013optimal}) an accelerated stochastic approximation
method and showed that the iteration complexity of this algorithm is bounded by
\[\cO\left(\sqrt{\frac{L_G}{\epsilon}} + \frac{\sigma_G^2}{\epsilon^2}\right).\]
More recently, Chen, Lan and Ouyang~\cite{chen2013optimal} presented an optimal stochastic accelerated primal-dual (APD) method 
with a better complexity bound than \eqref{eqnBestRate}
for solving $SVI(Z;G, H, J)$ with a linear operator $H$.

\subsection{Contribution of this paper}
Our contribution in this paper mainly consists of the following three aspects.
Firstly, we present the accelerated mirror-prox (AMP) method that computes a weak solution of $VI(Z; G, H, J)$ after incorporating a multi-step acceleration scheme into the mirror-prox method in \cite{nemirovski2005prox}. By utilizing the smoothness of $G(\cdot)$, we can significantly improve the iteration complexity from \eqref{eqnRateVI} to \eqref{eqnOptRate}, while the iteration cost of AMP is comparable to that of the mirror-prox method. Therefore, AMP can solve VI problems efficiently with big Lipschitz constant $L_G$. To the best of our knowledge, this is the first time in the literature that such an optimal iteration complexity bound has been obtained for general
Lipschitz continuous (rather than linear) operator $H$.

Secondly, we develop a stochastic counterpart of AMP, namely stochastic AMP, to compute a weak solution of $SVI(Z; G, H, J)$, and
demonstrate that its iteration complexity is bounded by \eqref{eqnOptRateS} and, similarly to the stochastic mirror-prox method,
each iteration of this algorithm requires the calls to ${\cal SO}_G$ and ${\cal SO}_H$.
Therefore, this algorithm improves the best-known complexity bounds for stochastic VI in terms of
the dependence on the Lipschitz constant $L_G$. To the best of our knowledge, this is the first time that such an optimal iteration complexity bound 
has been developed for $SVI(Z; G, H, J)$ for general Lipschitz continuous (rather than linear) operator $H$.  In addition, we investigate the stochastic VI method in more details, e.g., by developing the large-deviation results associated with the convergence of stochastic AMP.

Finally, for both deterministic and stochastic VI, we demonstrate that the AMP can deal with the case when $Z$ is unbounded,
as long as a strong solution to problem \eqref{eqnSVI} exists. We incorporate into AMP
the termination criterion employed by Monteiro and Svaiter~\cite{monteiro2010complexity, monteiro2011complexity} for solving variational and hemivariational inequalities
posed as monotone inclusion problem. In both deterministic and stochastic cases, when $Z$ is unbounded, the iteration complexity of AMP will depend on the distance from the initial point to the set of strong solutions.

\subsection{Organization of the paper}
The paper is organized as follows. We propose the AMP algorithm and discuss the main convergence results for solving deterministic VI and stochastic VI in Sections \ref{secAMP} and \ref{secAMPS} respectively. To facilitate the readers, we present the proofs of the main convergence results in Section \ref{secProof}. Finally, we make some concluding remarks in Section \ref{secConclusion}.

\section{Accelerated prox-method for deterministic VI}
\label{secAMP}

We introduce in this section an accelerated mirror-prox (AMP) method that computes a weak solution of $VI(Z;G, H, J)$, and discuss its main convergence properties. 

Throughout this paper, we assume that the following \emph{prox-mapping} can be solved efficiently:
\begin{equation}
	\label{eqnProxMapping}
	P_z^J(\eta):=\argmin_{u\in Z}\langle\eta, u-z\rangle + V(z, u) + J(u).
\end{equation}
In \eqref{eqnProxMapping}, the function $V(\cdot,\cdot)$ is defined by
\begin{align}
	\label{eqnV}
	V(z,u):=\omega(u)-\omega(z)-\langle\nabla\omega(z),u-z\rangle,\ \forall u,z\in Z,
\end{align}
where $\omega(\cdot)$ is a strongly convex function with convexity parameter $\mu>0$. The function $V(\cdot,\cdot)$ is known as a \emph{prox-function}, or \emph{Bregman divergence} \cite{bregman1967relaxation} (see, e.g., \cite{nemirovski2005prox,ben-tal2005non,nesterov2005smooth,auslender2006interior} for the properties of prox-functions and prox-mappings and their applications in convex optimization). Using the aforementioned definition of the prox-mapping, we describe the AMP method in Algorithm \ref{algAMP}.
\begin{algorithm}
	\caption{The accelerated mirror-prox (AMP) method for solving a weak solution of $VI(Z;G,H,J)$}
	\label{algAMP}
	\begin{algorithmic}
		\STATE Choose $r_1\in Z$. Set $w_1=r_1$, $w_1^{ag}=r_1$.
		\STATE For $t=1,2,\ldots, N-1$, calculate
			\begin{align}
				w\md  & = (1-\alpha_t)w\ag + \alpha_t r_t,\label{eqnwmd}\\
				w_{t+1} & = P_{r_t}^{\gamma_tJ}\left(\gamma_tH(r_t)+\gamma_t\nabla G(w\md)\right),\label{eqnProxwmd}\\
				r_{t+1} & = P_{r_t}^{\gamma_tJ}\left(\gamma_tH(w_{t+1})+\gamma_t\nabla G(w\md)\right),\label{eqnProxrmd}\\
				w\ag[1]  & = (1-\alpha_t)w\ag + \alpha_t w_{t+1}.\label{eqnwag}
			\end{align}
		\STATE Output $w_{N+1}^{ag}$.
	\end{algorithmic}
\end{algorithm}

Observe that the AMP method differs from the mirror-prox method in that we introduced two new sequences, i.e., $\{w\md \}$ and $\{w\ag\}$ (here ``md'' stands for ``middle'', and ``ag'' stands for ``aggregated''). On the other hand, the mirror-prox method only had to compute the ergodic mean of the sequence $\{w_t\}$ as the output of the algorithm (similar to $\{w\ag\}$). If $\alpha_t\equiv 1$,  $G=0$ and $J=0$, then Algorithm \ref{algAMP} for solving $VI(Z;0, H, 0)$ is equivalent to the prox-method in \cite{nemirovski2005prox}. In addition, if the distance generating function $w(\cdot)=\|\cdot\|^2/2$, then iterations \eqref{eqnProxwmd} and \eqref{eqnProxrmd} becomes
\begin{align*}
	w\tn & = \argmin_{u\in Z}\langle\gamma_tH(r_t), u-r_t\rangle + \frac 12\|u-r_t\|^2,
	\\
	r\tn & = \argmin_{u\in Z}\langle\gamma_tH(w\tn), u-r_t\rangle + \frac 12\|u-r_t\|^2,
\end{align*}
which are exactly the iterates of the extragradient method in \cite{korpelevich1976extragradient}. On the other hand, if $H=0$, then the iterations \eqref{eqnProxwmd} and \eqref{eqnProxrmd} produces the same optimizer $w\tn = r\tn$, and Algorithm \ref{algAMP} is equivalent to a version of Nesterov's accelerated method for solving $\min_{u\in Z}G(u)+J(u)$ (see, for example, Algorithm 1 in  \cite{tseng2008accelerated}). Therefore, Algorithm \ref{algAMP} can be viewed as a hybrid algorithm of the mirror-prox method and the accelerated gradient method, which gives its name accelerated mirror-prox method.

In order to analyze the convergence of Algorithm \ref{algAMP}, we introduce a notion to characterize the weak solutions of $VI(Z;G,H,J)$. For all $\tilde u, u\in Z$, we define
\begin{align}
	\label{eqnQ}
	Q(\tilde u, u):=G(\tilde u) - G(u) + \langle H(u), \tilde u - u\rangle + J(\tilde u) - J(u).
\end{align}
Clearly, for $F$ defined in \eqref{eqnF}, we have $\langle F(u),\tilde{u}-u\rangle \le Q(\tilde{u}, u)$. Therefore, if $Q(\tilde u, u)\leq 0$ for all $u\in Z$, then $\tilde u$ is a weak solution of $VI(Z;G,H,J)$. Hence when $Z$ is bounded, it is natural to use the gap function
\begin{equation}
	\label{eqng0}
	g(\tilde{u}):= \sup_{u\in Z}Q(\tilde u, u)
\end{equation}
to evaluate the quality of a feasible solution $\tilde{u}\in Z$. However, if $Z$ is unbounded, then $g(\tilde{z})$ may not be well-defined, even when $\tilde{z}\in Z$ is a nearly optimal solution. Therefore, we need to employ a slightly modified gap function in order to measure the quality of candidate solutions when $Z$ is unbounded. In the sequel, we will consider the cases of bounded and unbounded $Z$ separately.

\vgap

Theorem \ref{thmAMPRateB} below describes the convergence property of Algorithm \ref{algAMP} when $Z$ is bounded. It should be noted that the following quantity will be used throughout the convergence analysis of this paper:
		\begin{equation}
			\label{eqnGamma}
			\Gamma_t = \left\{\begin{aligned}
				&1, & \text{ when } t=1\\
				&(1-\alpha_t)\Gamma_{t-1}, & \text{ when } t>1,\\
			\end{aligned} \right.
		\end{equation}

\begin{thm}
	\label{thmAMPRateB}
	Suppose that 
	\begin{equation}
		\label{eqnVBounded}
		\ds\sup_{z_1, z_2\in Z}V(z_1, z_2) \le \Omega_Z^2. 
	\end{equation}
	If the parameters $\{\alpha_t\}$ and $\{\gamma_t\}$ in Algorithm \ref{algAMP} are chosen such that $\alpha_1=1$, and 
	\begin{equation}
		\label{eqnCondAlphaGammaBD}
		0\leq \alpha_t<1, \ {\mu}-{L_G\alpha_t\gamma_t}- \frac{L_H^2\gamma_t^2}{\mu}\geq 0, \textrm{ and }\frac{\alpha_t}{\Gamma_t\gamma_t} \leq \frac{\alpha_{t+1}}{\Gamma_{t+1}\gamma_{t+1}},\ \forall t\ge 1,
	\end{equation}
	where $\{\Gamma_t\}$ is defined by \eqref{eqnGamma}. Then,
	\begin{equation}
		\label{eqngBound}
		\ds g(w\ag[1])\leq \frac{\alpha_t}{\gamma_t}\Omega_Z^2.
	\end{equation}	
\end{thm}

There are various options for choosing the parameters $\{\alpha_t\}$ and $\{\gamma_t\}$ that satisfy \eqref{eqnCondAlphaGammaBD}. In the following corollary, we give one example of such parameter settings.

\begin{cor}
	\label{corAMPRateB}
	Suppose that \eqref{eqnVBounded} holds. If the parameters $\{\alpha_t\}$ and $\{\gamma_t\}$ in AMP are set to
	\begin{align}
		\label{eqnAlphaGammaBounded}
		\alpha_t = \frac 2{t+1}\ \text{and }
		\gamma_t  = \frac{\mu t}{2(L_G+L_Ht)},
	\end{align}
	then for all $u\in Z$, 
		\begin{align}
			\label{eqnAMPRateB}
			&	Q(w\ag[1], u) 
			\leq  \left(\frac{4L_G}{\mu t(t+1)} + \frac{4L_H}{\mu t} \right)\Omega_Z^2,
		\end{align}		
	where $\Omega_Z$ is defined in \eqref{eqnVBounded}.
	\begin{proof}
		Clearly, $\ds\Gamma_t=\frac{2}{t(t+1)}$ satisfies \eqref{eqnGamma}, and
		\[		\frac{\alpha_t}{\Gamma_t\gamma_t} = \frac{2}{\mu}(L_G + L_Ht)\leq \frac{\alpha_{t+1}}{\Gamma_{t+1}\gamma_{t+1}}. \]
		Moreover,
		\begin{align*}
			& {\mu}-{L_G\alpha_t\gamma_t}- \frac{L_H^2\gamma_t^2}{\mu} = \mu - \frac{\mu L_G}{L_G + L_Ht}\cdot \frac t {t+1} - \frac{\mu L_H^2t^2}{4(L_G + L_Ht)^2}
			\geq \mu - \frac{\mu L_G}{L_G + L_Ht} - \frac{\mu L_Ht}{L_G + L_Ht} = 0.
		\end{align*}
		Thus \eqref{eqnCondAlphaGammaBD} holds. Hence, by applying \eqref{eqngBound} in Theorem \ref{thmAMPRateB} with the parameter setting in \eqref{eqnAlphaGammaBounded} and using \eqref{eqnVBounded}, we obtain \eqref{eqnAMPRateB}.
	\end{proof}
\end{cor}
\vgap

Clearly, in view of \eqref{eqnAMPRateB}, when the parameters are chosen according to \eqref{eqnAlphaGammaBounded}, 
the number of iterations performed by the AMP method to find an $\epsilon$-solution of \eqref{eqnProblem},
i.e., a point $\bar w \in Z$ s.t. $g(\bar w) \le \epsilon$, can be bounded
by 
\[
{\cal O} \left( \sqrt{\frac{L_G}{\epsilon}} + \frac{L_H}{\epsilon} \right).
\]
This bound significantly improves the best-known so-far complexity for solving problem \eqref{eqnProblem} (see \eqref{eqnOptRate})
in terms of their dependence on the Lipschitz constant $L_G$. Moreover, it should be noted that the parameter setting in \eqref{eqnAlphaGammaBounded} is independent of $\Omega_Z$, i.e., the AMP method achieves the above optimal iteration-complexity without requiring any information on the diameter of $Z$.

\vgap 

Now, we consider the case when $Z$ is unbounded. To study the convergence properties of AMP in this case, we use a perturbation-based termination criterion recently employed by Monteiro and Svaiter \cite{monteiro2010complexity, monteiro2011complexity}, which is based on the enlargement of a maximal monotone operator first introduced in \cite{burachik1997enlargement}. More specifically, we say that the pair $(\tilde v, \tilde u)\in \cE\times Z$ is a $(\rho, \varepsilon)$-approximate solution of $VI(Z;G,H,J)$ if $\|\tilde v\|\leq \rho$ and $\tilde g(\tilde u, \tilde v)\leq \varepsilon$, where the gap function $\tilde g(\cdot, \cdot)$ is defined by
\begin{equation}
	\label{eqngt}
	\tilde g(\tilde u, \tilde v):=\sup_{u\in Z}Q(\tilde u, u) - \langle\tilde v, \tilde u - u\rangle.
\end{equation}
We call $\tilde v$ the perturbation vector associated with $\tilde u$. One advantage of employing this termination criterion is that the convergence analysis does not depend on the boundedness of $Z$. 

Theorem \ref{thmAMPRateUB} below describes the convergence properties of AMP for solving deterministic VIs with unbounded feasible sets, under the assumption that a strong solution of \eqref{eqnProblem} exists. It should be noted that this assumption does not limit too much the applicability of the AMP method. For example, when $J(\cdot)=0$, any weak solution to $VI(Z;F)$ is also a strong solution.

\begin{thm}
	\label{thmAMPRateUB}
	Suppose that $V(r, z):= \|z - r\|^2/2$ for any $r, z\in Z$. Also assume that the parameters $\{\alpha_t\}$ and $\{\gamma_t\}$ in Algorithm \ref{algAMP} are chosen such that $\alpha_1=1$, and for all  $t> 1$,
	\begin{equation}
		\label{eqnCondAlphaGammaUB}
		0\leq \alpha_t<1, \ {L_G\alpha_t\gamma_t}+{L_H^2\gamma_t^2}\leq c^2 \text{ for some }c<1, \text{ and }\frac{\alpha_t}{\Gamma_t\gamma_t} = \frac{\alpha_{t+1}}{\Gamma_{t+1}\gamma_{t+1}},
	\end{equation}	
	where $\Gamma_t$ is defined in \eqref{eqnGamma}. 
	Then for all $t\ge 1$ there exist $v\tn\in \cE$ and $\varepsilon\tn\geq 0$ such that $\tilde g(w\ag[1], v\tn) \leq\varepsilon\tn$. Moreover, we have
	\begin{align}
		\label{eqnvepsThm}
		\|v\tn\| \leq \frac{2\alpha_tD}{\gamma_t}\text{ and }
		\varepsilon\tn\leq \frac{3\alpha_t(1+\theta_t)D^2}{\gamma_t}.
	\end{align}
	where 
	\begin{align}
		\label{eqndtheta}
		D: = \|r_1 - u^*\|,\ 
		\theta_t:=\frac{\Gamma_t}{2(1-c^2)}\max_{i=1,\ldots,t}\frac{\alpha_i}{\gamma_i}.
	\end{align}
	and $u^*$ is a strong solution of $VI(Z;G,H,J)$.
\end{thm}

\vgap 

Below we provide a specific setting of parameters $\{\alpha_t\}$ and $\{\gamma_t\}$ that satisfies condition \eqref{eqnCondAlphaGammaUB}.

\begin{cor}
	\label{corAMPRateUB}
	Suppose that $V(r, z):= \|z - r\|^2/2$ for any $r, z\in Z$ and $L_H>0$. In Algorithm \ref{algAMP}, if $N\ge 2$ is given and the parameters $\{\alpha_t\}$ and $\{\gamma_t\}$ are set to
	\begin{align}
		\label{eqnAlphaGammaUB}
		\alpha_t = \frac 2{t+1}\ \text{ and }
		\gamma_t  = \frac{t}{3(L_G+L_HN)},
	\end{align}
	then there exists $v_N\in \cE$ such that $\tilde g(w^{ag}_N, v_N)\leq \varepsilon_N$, 
	\begin{align}
		\label{eqnvepsCor}
		\|v_N\| \leq \left[\frac{12L_G}{N(N-1)} + \frac{12L_H}{N-1}\right]D,
		\text{ and }
		\varepsilon_N \leq  \left[\frac{45L_G}{N(N-1)} + \frac{45L_H}{N-1}\right]D^2,
	\end{align}
	where $u^*$ is a strong solution of $VI(Z;F)$ and $D$ is defined in \eqref{eqndtheta}.
	\begin{proof}
		Clearly, $\ds\Gamma_t=\frac{2}{t(t+1)}$ satisfies \eqref{eqnGamma}, and 
		\begin{align*}
			L_G\alpha_t\gamma_t + L_H^2\gamma_t^2 = &\ \frac{2L_Gt}{3(L_G+L_HN)(t+1)} + \frac{L_H^2t^2}{9(L_G + L_HN)^2} \leq \frac{2L_G}{3(L_G+L_HN)} + \frac{L_HN}{3(L_G + L_HN)}
			\\
			=& \ \frac{2L_G + L_HN}{3(L_G + L_HN)}\le \frac{2}{3} =: c^2.
		\end{align*}
		We can see that $c<1$, and $\ds\frac{1}{1-c^2}=3$. Moreover, when $N\ge 2$,
		\begin{align*}
			\theta_{N-1} = {\frac{\Gamma_{N-1}}{2(1-c^2)}}\max_{1\leq i\leq {N-1}}\{\frac{\alpha_i}{\Gamma_i}\} = {\frac{1}{(1-c^2)N(N-1)}}\max_{1\leq i\leq {N-1}}i = \frac{3}{N} \le \frac{3}{2}.
		\end{align*}
		We conclude \eqref{eqnvepsCor} by substituting the values of $\alpha_{N-1}$,  $\gamma_{N-1}$ and $\theta_{N-1}$ to \eqref{eqnvepsThm}.
	\end{proof}
\end{cor}
\vgap

Several remarks are in place for the results obtained in Theorem \ref{thmAMPRateUB} and Corollary \ref{corAMPRateUB}. Firstly, although the existence of a strong solution $u^*$ is assumed, no information on either $u^*$ or $D$ is needed for choosing parameters $\alpha_t$ and $\gamma_t$, as shown in \eqref{eqnAlphaGammaUB} of Corollary \ref{corAMPRateUB}. Secondly, both residuals $\|v_N\|$ and $\varepsilon_N$ in \eqref{eqnvepsCor} converge to $0$ at the same rate (up to a constant $15D/4$). Finally, it is only for simplicity that we assume that $V(r,z)=\|z-r\|^2/2$; Similar results can be achieved under assumptions that $\nabla\omega$ is Lipschitz continuous and that $\sqrt{V(\cdot, \cdot)}$ is a metric.

\section{Accelerated prox-method for stochastic VI}
\label{secAMPS}

In this section, we focus on the $SVI(Z;F)$, and demonstrate that the stochastic counterpart of Algorithm \ref{algAMP} can achieve the optimal rate of convergence in \eqref{eqnOptRateS}.

The stochastic AMP is obtained by replacing the operators $H(r_t), H(w\tn)$ and $\nabla G(x\md)$ in Algorithm \ref{algAMP} by their stochastic counterparts $\cH(r_t, \zeta_{2t-1})$, $\cH(w_{t+1}, \zeta_{2t})$ and $\cG(w\md, \xi_t)$ respectively, by calling the stochastic oracles ${\cal SO}_G$ and ${\cal SO}_H$. This algorithm is formally described in Algorithm \ref{algAMPS}.
\begin{algorithm}[!h]
	\caption{The accelerated mirror-prox (AMP) method for solving a weak solution of $SVI(Z;G,H,J)$}
	\label{algAMPS}
	\begin{algorithmic}
		\STATE Modify \eqref{eqnProxwmd} and \eqref{eqnProxrmd} in Algorithm \ref{algAMP} to
			\begin{align}
				w_{t+1} & = P_{r_t}^{\gamma_tJ}\left(\gamma_t\cH(r_t, \zeta_{2t-1})+\gamma_t \cG(w\md, \xi_t)\right),\label{eqnProxwmdS}\\
				r_{t+1} & = P_{r_t}^{\gamma_tJ}\left(\gamma_t\cH(w_{t+1}, \zeta_{2t})+\gamma_t\cG(w\md, \xi_t)\right),\label{eqnProxrmdS}
			\end{align}
	\end{algorithmic}
\end{algorithm}

It is interesting to note that for any $t$, there are two calls of ${\cal SO}_H$ but just one call of ${\cal SO}_G$. However, if we assume that $J=0$ and use the stochastic mirror-prox method in \cite{juditsky2011solving} to solve $SVI(Z;G,H,0)$, for any $t$ there would be two calls of ${\cal SO}_H$ and two calls of ${\cal SO}_G$. Therefore, the cost per iteration of AMP is less than that of the stochastic mirror-prox method.

Similarly to Section \ref{secAMP}, we use the gap function $g(\cdot)$ for the case when $Z$ is bounded, and use the modified gap function $\tilde g(\cdot,\cdot)$ for the case when $Z$ is unbounded. For both cases we establish the rate of convergence of the gap functions in terms of their expectation, i.e., the ``average'' rate of convergence over many runs of the algorithm. Furthermore, we demonstrate that if $Z$ is bounded, then we can also establish the rate of convergence of $g(\cdot)$ in the probability sense, under the following ``light-tail'' assumption:
\begin{enumerate}
	\setcounter{enumi}{\theenuAssumptions}
	\renewcommand{\theenumi}{\textbf{A\arabic{enumi}}}
	\item 
	\label{itmLT}
	For any $i$-th call on oracles ${\cal SO}_H$ and ${\cal SO}_H$ with any input $u\in Z$, 
	$$\E[\exp\{\|\nabla G(u) - \cG(u, \xi_i)\|_*^2/\sigma_G^2 \}]\leq \exp\{1\},\text{ and }\E[\exp\{\|H(u) - \cH(u, \zeta_i)\|_*^2/\sigma_H^2 \}]\leq \exp\{1\}. $$
	\setcounter{enuAssumptions}{\value{enumi}}
\end{enumerate}
It should be noted that Assumption \ref{itmLT} implies Assumption \ref{itmVB} by Jensen's inequality.

\vgap 

The following theorem shows the convergence property of Algorithm \ref{algAMPS} when $Z$ is bounded.
\begin{thm}
	\label{thmAMPRateBS}
	Suppose that \eqref{eqnVBounded} holds. Also assume that the parameters $\{\alpha_t\}$ and $\{\gamma_t\}$ in Algorithm \ref{algAMPS} satisfy $\alpha_1=1$ and
	\begin{align}
		\label{eqnCondAlphaGammaIncS}
		q{\mu}-{L_G\alpha_t\gamma_t}- \frac{3L_H^2\gamma_t^2}{\mu}\geq 0 \text{ for some }q\in(0,1) ,
		\text{ and } \frac{\alpha_t}{\Gamma_t\gamma_t} \leq \frac{\alpha_{t+1}}{\Gamma_{t+1}\gamma_{t+1}},	\ \forall t\geq 1,
	\end{align}
	where $\Gamma_t$ is defined in \eqref{eqnGamma}. 
	Then,
	\begin{enumerate}
	\renewcommand{\theenumi}{(\alph{enumi})}
	\renewcommand{\labelenumi}{\theenumi\!\!}
	\item Under Assumption \ref{itmVB}, for all $t\ge 1$,
	\begin{align}
		\label{eqnQBoundBS}
		&\E\left[g(w\ag[1])\right] \leq \cQ_0(t),
		\\
		\label{eqnQ0}
		\textrm{where }&\cQ_0(t):=\frac{2\alpha_t}{\gamma_t}\Omega_Z^2 + \left[ 4\sigma_H^2 + \left(1 + \frac{1}{2(1-q)} \right)\sigma_G^2\right]\Gamma_t
				\sum_{i=1}^{t}\frac{\alpha_i\gamma_i}{\mu\Gamma_i}.
	\end{align}		
	\item Under Assumption \ref{itmLT}, for all $\lambda>0$ and $t\ge 1$,
	\begin{align}
		\label{eqnProb}
		& Prob\{g(w\ag[1])>\cQ_0(t) + \lambda\cQ_1(t) \}\leq 2\exp\{-\lambda^2/3 \} + 3\exp\{-\lambda \}
		\\
		\label{eqnQ1}
		\textrm{where }&\cQ_1(t):= \Gamma_t(\sigma_G+\sigma_H) \Omega_Z\sqrt{\frac{2}{\mu}\sum_{i=1}^t\left(\frac{\alpha_i}{\Gamma_i} \right)^2 } + \left[ 4\sigma_H^2 + \left(1 + \frac{1}{2(1-q)} \right)\sigma_G^2\right]\Gamma_t\sum_{i=1}^{t}\frac{\alpha_i\gamma_i}{\mu\Gamma_i}.
	\end{align}
	\end{enumerate}
\end{thm}

We present below a specific parameter setting of $\{\alpha_t\}$ and $\{\gamma_t\}$ that satisfies \eqref{eqnCondAlphaGammaIncS}.
\begin{cor}
	\label{corStepS}
	Suppose that \eqref{eqnVBounded} holds. If the stepsizes $\{\alpha_t\}$ and $\{\gamma_t\}$ in Algorithm \ref{algAMPS} are set to:
	\begin{align}
		\label{eqnAlphaGammaBDS}
		\alpha_t = \frac 2{t+1}\text{ and } 
		\gamma_t = \frac{\mu t}{4L_G+ 3L_Ht +  \sigma (t+1)\sqrt{\mu t}/(\sqrt{2}\Omega_Z)},
	\end{align}
	where
	$
		\sigma:= \sqrt{\sigma_H^2 + \sigma_G^2},
	$
	and $\Omega_Z$ is defined in \eqref{eqnVBounded}. Then under Assumption \ref{itmVB},
		\begin{align}
			\label{eqnC0}
			&\ \E\left[g(w\ag[1])\right] 
			\leq \frac{16L_G \Omega_Z^2}{\mu t(t+1)} + \frac{12L_H \Omega_Z^2}{\mu(t+1)} + \frac{7(\sigma_G+\sigma_H) \Omega_Z}{\sqrt{\mu(t-1)}} =:\cC_0(t).
		\end{align}
	Furthermore, under Assumption \ref{itmLT}, 
	\[Prob\{g(w\ag[1])>\cC_0(t) + \lambda\cC_1(t) \}\leq 2\exp\{-\lambda^2/3 \} + 3\exp\{-\lambda \}, \ \forall\lambda>0,\]
	where
	\begin{align}
		\label{eqnC1}
		\cC_1(t) := \frac{6(\sigma_G + \sigma_H)\Omega_Z}{\sqrt{\mu(t-1)}}.
	\end{align}
	
	\begin{proof}
		It is easy to check that $\ds\Gamma_t=\frac{2}{t(t+1)}$ and $\ds\frac{\alpha_t}{\Gamma_t\gamma_t} \leq \frac{\alpha_{t+1}}{\Gamma_{t+1}\gamma_{t+1}}$. In addition, in view of \eqref{eqnAlphaGammaBDS}, we have $\gamma_t\leq \mu t/(4L_G)$ and $\gamma_t^2 \leq (\mu^2)/(9L_H^2)$, which implies
		\begin{align*}
			\frac{5\mu}{6} - L_G\alpha_t\gamma_t - \frac{3L_H^2\gamma_t^2}{\mu} \geq \frac{5\mu}{6} - \frac{\mu t}{4}\cdot\frac{2}{t+1} - \frac{\mu}{3} \geq 0.
		\end{align*}
		Therefore the first relation in \eqref{eqnCondAlphaGammaIncS} holds with constant $q=5/6$. In view of Theorem \ref{thmAMPRateBS}, it now suffices to show that $\cQ_0(t)\leq \cC_0(t)$ and $\cQ_1(t)\leq \cC_1(t)$. Observe that $\alpha_t/\Gamma_t=t$, and $\gamma_t\leq (\sqrt{2\mu} \Omega_Zt)/(\sigma t^{3/2})$, thus using the fact that $\sum_{i=1}^t\sqrt{i}\le \int_{0}^{t+1}\sqrt{t}dt = \frac{2}{3}(t+1)^{3/2}$, we obtain
		\begin{align*}
			&\ \sum_{i=1}^{t}\frac{\alpha_i\gamma_i}{\Gamma_i} \leq \sum_{i=1}^t \frac{\sqrt{2\mu} \Omega_Z i^2}{\sigma i^{3/2}} = \frac{\sqrt{2\mu} \Omega_Z}{\sigma}\sum_{i=1}^t\sqrt{i}\leq \frac{\sqrt{2\mu} \Omega_Z}{3\sigma}(t+1)^{3/2}.
		\end{align*}
		Using the above inequality, \eqref{eqnVBounded}, \eqref{eqnQ0}, \eqref{eqnQ1}, \eqref{eqnAlphaGammaBDS}, and the fact that $\sqrt{t+1}/t\leq 1/\sqrt{t-1}$ and $\sum_{i=1}^t i^2\leq t(t+1)^2/3$, we have
		\begin{align*}
			\cQ_0(t)= &\ \frac{4\Omega_Z^2}{\mu t(t+1)}\left(4L_G+ 3L_Ht + \sigma (t+1)\sqrt{2\mu t}/\Omega_Z \right) + \frac{4\sigma^2}{ \mu t(t+1)}\sum_{i=1}^{t}\frac{\alpha_i\gamma_i}{\Gamma_i}
			\\
			\leq &\ \frac{16L_G \Omega_Z^2}{\mu t(t+1)} + \frac{12L_H \Omega_Z^2}{\mu(t+1)} + \frac{2\sqrt{2}\sigma \Omega_Z}{\sqrt{\mu t}} + \frac{8\sigma \Omega_Z\sqrt{2(t+1)}}{3\sqrt{\mu}t}
			\leq \cC_0(t),
		\end{align*}
		and
		\begin{align*}
			\cQ_1(t) = &\ \frac{2(\sigma_G + \sigma_H)}{t(t+1)} \Omega_Z\sqrt{\frac{2}{\mu}\sum_{i=1}^ti^2} + \frac{8\sigma^2}{ \mu t(t+1)}\sum_{i=1}^{t}\frac{\alpha_i\gamma_i}{\Gamma_i}
			\leq \frac{2\sqrt{2}(\sigma_G + \sigma_H)\Omega_Z}{\sqrt{3\mu t}} + \frac{8\sigma \Omega_Z\sqrt{2(t+1)}}{3\sqrt{\mu}t}
			\leq \cC_1(t).
		\end{align*}
	\end{proof}
\end{cor}

\vgap

In view of \eqref{eqnOptRateS}, \eqref{eqnC0} and \eqref{eqnC1}, we can clearly see that the stochastic AMP method achieves the optimal iteration complexity for solving the SVI problem. More specifically, this algorithm allows $L_G$ to be as large as $\cO(t^{3/2})$ without significantly affecting its convergence properties.

\vgap 

In the following theorem, we demonstrate the convergence properties of Algorithm \ref{algAMPS} for solving the stochastic problem $SVI(Z;G,H,J)$ when $Z$ is unbounded. It seems that this case has not been well-studied previously in the literature.

\begin{thm}
	\label{thmAMPRateUBS}
	Suppose that $V(r, z):= \|z - r\|^2/2$ for any $r\in Z$ and $z\in Z$. If the parameters $\{\alpha_t\}$ and $\{\gamma_t\}$ in Algorithm \ref{algAMP} are chosen such that $\alpha_1=1$, and for all $t> 1$,
	\begin{equation}
		\label{eqnCondAlphaGammaUBS}
		0\leq \alpha_t<1, \ {L_G\alpha_t\gamma_t}+{3L_H^2\gamma_t^2}\leq c^2 < q \text{ for some }c,q\in(0,1), \text{ and }\frac{\alpha_t}{\Gamma_t\gamma_t} = \frac{\alpha_{t+1}}{\Gamma_{t+1}\gamma_{t+1}},
	\end{equation}	
	where $\Gamma_t$ is defined in \eqref{eqnGamma}. 
	Then for all $t\ge 1$ there exists a perturbation vector $v\tn$ and a residual $\varepsilon\tn\geq 0$ such that $\tilde g(w\ag[1], v\tn) \leq\varepsilon\tn$. Moreover, for all $t\ge 1$, we have
	\begin{align}
		\label{eqnEv}
		\E[\|v\tn\|] \leq&\ \frac{\alpha_t}{\gamma_t}(2D + 2\sqrt{D^2 + C_t^2}),
		\\ 
		\label{eqnEeps}
		\E[\varepsilon\tn]\leq&\ \frac{\alpha_t}{\gamma_t}\left[(3+6\theta)D^2 + (1+6\theta)C_t^2\right] + \frac{18\alpha_t^2\sigma_H^2}{\gamma_t^2}\sum_{i=1}^t\gamma_i^3,
	\end{align}
	where $u^*$ is a strong solution of $VI(Z;G,H,J)$,
		\begin{align}
		\label{eqnCtheta}
		\theta = \max\{1, \frac{c^2}{q-c^2} \}\text{ and } C_t = \sqrt{\left[ 4\sigma_H^2 + \left(1 + \frac{1}{2(1-q)} \right)\sigma_G^2\right]\sum_{i=1}^{t}\gamma_i^2}.
		\end{align}
\end{thm}

Below we give an example of parameters $\alpha_t$ and $\gamma_t$ that satisfies \eqref{eqnCondAlphaGammaUBS}.
\begin{cor}
	\label{corStepUBS}
	Suppose that there exists a strong solution of \eqref{eqnProblem}. If the maximum number of iterations $N$ is given, and the stepsizes $\{\alpha_t\}$ and $\{\gamma_t\}$ in Algorithm \ref{algAMPS} are set to
	\begin{align}
		\label{eqnStepUBS}
		\alpha_t = \frac 2{t+1}\text{ and }
		\gamma_t = \frac{t}{5L_G + 3L_HN + \sigma N\sqrt{N-1}/\tilde{D}},
	\end{align}
	where $\sigma$ is defined in Corollary \ref{corStepS}, then there exists $v_N\in \cE$ and $\varepsilon_N>0$, such that $\tilde{g}(w\ag[][N],v_{N})\le\varepsilon_N$, 
	\begin{align}
		\label{eqnEvUB}
		\E[\|v_N\|]\leq&\  \frac{40L_GD}{N(N-1)} + \frac{24L_HD}{N-1} + \frac{\sigma(8D/\tilde{D} + 5)}{\sqrt{N-1}},
	\end{align}
	and
	\begin{align}
		\label{eqnEepsUB}
		\E[\varepsilon_N]\le &\ \frac{90L_GD^2}{N(N-1)} + \frac{54L_HD^2}{N-1} + \frac{\sigma D}{\sqrt{N-1}}\left(\frac{18D}{\tilde{D}}+ (19+\frac{18}{N})\frac{\tilde{D}}{D}\right).
	\end{align}
	\begin{proof}
		Clearly, $\Gamma_t =\ds\frac{2}{t(t+1)}$ satisfies \eqref{eqnGamma}. Moreover, in view of \eqref{eqnStepUBS}, we have
		\begin{align*}
			L_G\alpha_t\gamma_t + 3L_H^2\gamma_t^2 \le \frac{2L_G}{5L_G+3L_HN} + \frac{3L_H^2N^2}{(5L_G+3L_HN)^2} = \frac{10L_G^2 + 6L_GL_HN + 3L_H^2N^2}{(5L_G+3L_HN)^2} < \frac{5}{12} < \frac{5}{6},
		\end{align*}
		which implies that  \eqref{eqnCondAlphaGammaUBS} is satisfied with $c^2=5/12$ and $q=5/6$. Observing from \eqref{eqnStepUBS} that $\gamma_t = t\gamma_1$, setting $t=N-1$ in \eqref{eqnCtheta} and \eqref{eqnStepUBS}, we obtain
		\begin{align}
			\label{eqnCN}
			\frac{\alpha_{N-1}}{\gamma_{N-1}} = \frac{2}{\gamma_1N(N-1)}\text{ and } C_{N-1}^2 = 4\sigma^2\sumt[N-1] \gamma_1^2 i^2 \le \frac{4\sigma^2\gamma_1^2N^2(N-1)}{3}.
		\end{align}
		Applying \eqref{eqnCN} to \eqref{eqnEv} we have
		\begin{align*}
			\E[\|v_N\|] \le&\ \frac{2}{\gamma_1N(N-1)}(4D + 2C_{N-1}) \le \frac{8D}{\gamma_1N(N-1)} + \frac{8\sigma}{\sqrt{3(N-1)}}
			\\
			\le &\ \frac{40L_GD}{N(N-1)} + \frac{24L_HD}{N-1} + \frac{\sigma(8D/\tilde{D} + 5)}{\sqrt{N-1}}.
		\end{align*}
		In addition, using  \eqref{eqnEeps}, \eqref{eqnCN}, and the facts that $\theta=1$ in \eqref{eqnCtheta} and $\sumt[N-1]\gamma_i^3 = N^2(N-1)^2/4$, we have
		\begin{align*}
			&\ \E[\varepsilon_{N-1}] \le \frac{2}{\gamma_1N(N-1)}(9D^2 + 7C_{N-1}^2) + \frac{72\sigma_H^2}{\gamma_1^2N^2(N-1)^2}\cdot\frac{\gamma_1^3N^2(N-1)^2}{4}
			\\
			\le&\ \frac{18D^2}{\gamma_1 N(N-1)} + \frac{56\sigma\tilde{D}}{3\sqrt{N-1}} + \frac{18\sigma_H^2\tilde{D}}{\sigma N\sqrt{N-1}}
			\le\ \frac{90L_GD^2}{N(N-1)} + \frac{54L_HD^2}{N-1} + \frac{\sigma D}{\sqrt{N-1}}\left(\frac{18D}{\tilde{D}}+ (19+\frac{18}{N})\frac{\tilde{D}}{D}\right).
		\end{align*}
	\end{proof}
\end{cor}

\vgap

Observe that we need to choose a parameter $\tilde D$ for the stochastic unbounded case, which is not required for the deterministic case (see Corollary \ref{corAMPRateUB}). Note  that the value of $D$ will be very difficult to estimate for the unbounded case
and hence one often has to resort to a suboptimal selection for $\tilde D$.
For example, if $\tilde D = 1$, then the RHS of \eqref{eqnEvUB} and \eqref{eqnEepsUB}
will become ${\cal O}(L_G D/N^2 + L_HD/N + \sigma D / \sqrt{N})$
and ${\cal O}(L_GD^2/N^2 + L_HD^2/N + \sigma D^2 / \sqrt{N})$, respectively.

\section{Convergence analysis}
\label{secProof}
In this section, we focus on proving the main convergence results in Section \ref{secAMP} and \ref{secAMPS}, namely, Theorems \ref{thmAMPRateB}, \ref{thmAMPRateUB}, \ref{thmAMPRateBS} and \ref{thmAMPRateUBS}.

\subsection{Convergence analysis for deterministic AMP}
\label{secProofD}
In this section, we prove Theorems \ref{thmAMPRateB} and \ref{thmAMPRateUB} in Section \ref{secAMP}, which state the main convergence properties of Algorithm \ref{algAMP} for solving the deterministic problem $VI(Z;G,H,J)$.

To prove the convergence of the deterministic AMP, first we present some technical results. Propositions \ref{propProxMap} and \ref{propPRecursion} describe some important properties of the prox-mapping $P_r^J(\eta)$ in iterations \eqref{eqnProxwmd} and \eqref{eqnProxrmd} of Algorithm \ref{algAMP}. Proposition \ref{propSimplifiedQ} provides a recursion property of function $Q(\cdot,\cdot)$ defined in \eqref{eqnQ}. With the help of Propositions \ref{propProxMap}, \ref{propPRecursion} and \ref{propSimplifiedQ}, we can estimate a bound on $Q(\cdot,\cdot)$ in Lemma \ref{lemQBoundGeneral}.

\begin{pro}\label{propProxMap}
	For all $r, \zeta\in \cE$, if $w=P_r^J(\zeta)$, then for all $u\in Z$, we have
	\begin{equation}
		\label{eqnPOptCond}
		\langle\zeta, w-u\rangle + J(w) - J(u)\leq V(r, u)-V(r, w)-V(w,u).
	\end{equation}
	\begin{proof}
	 	See Lemma 2 in \cite{ghadimi2012optimal} for the proof.
	\end{proof}
\end{pro}

\vgap 

The following proposition is a slight extension of Lemma 6.3 in \cite{juditsky2011solving}. In particular, when $J(\cdot)=0$, we can obtain  \eqref{eqnPRecursion} and \eqref{eqnVvr} directly by applying \eqref{eqnLM} to (6.8) in \cite{juditsky2011solving}, and the results when $J(\cdot)\not\equiv 0$ can be easily constructed from the proof of Lemma 6.3 in \cite{juditsky2011solving}. We provide the proof here only for the integrity of this proposition.
\begin{pro}
	\label{propPRecursion}
	Given $r, w, y \in Z$ and $\eta, \vartheta\in \cE$ that satisfies
	\begin{align}
		\label{eqnProx1}
		w & = P_r^J(\eta),\\
		\label{eqnProx2}
		y & = P_r^J(\vartheta),
	\end{align}
	and
	\begin{align}
		\label{eqnLM}
		\|\vartheta - \eta\|^2_* \leq L^2\|w - r\|^2 + M^2,
	\end{align}	
	then for all $u\in Z$ we have
			\begin{align}
				\label{eqnPRecursion}
				\begin{aligned}
					&\ \langle \vartheta, w-u\rangle + J(w)-J(u)
					\leq V(r, u)-V(y, u)-\left(\frac{\mu}2- \frac{L^2}{2\mu}\right)\|r-w\|^2 + \frac{M^2}{2\mu},
				\end{aligned}
			\end{align}
			and
			\begin{align}
				\label{eqnVvr}
				V(y, w) \leq \frac{L^2}{\mu^2}V(r, w) + \frac{M^2}{2\mu}.
			\end{align}

	\begin{proof}
		Applying Proposition \ref{propProxMap} to \eqref{eqnProx1} and \eqref{eqnProx2}, for all $u\in Z$ we have
		\begin{align}
			\label{eqnMP2u}
			\langle \eta, w-u\rangle + J(w)-J(u) & \leq V(r, u)-V(r, w)-V(w, u),
			\\
			\label{eqnMP1}
			\langle \vartheta, y-u\rangle + J(y)-J(u) & \leq V(r, u)-V(r, y) - V(y, u),	
		\end{align}
		Specifically, letting $u=y$ in \eqref{eqnMP2u} we have
		\begin{align}\label{eqnMP2}
			\langle \eta, w-y\rangle + J(w)-J(y)\leq V(r, y)-V(r, w)-V(w, y).
		\end{align}
		Adding inequalities (\ref{eqnMP1}) and (\ref{eqnMP2}), then
		\begin{align*}
			\begin{aligned}
				&\ \langle\vartheta, y-u\rangle + \langle\eta, w-y\rangle + J(w)-J(u)
				\leq V(r, u)-V(y, u)-V(r,w)-V(w, y),
			\end{aligned}
		\end{align*}
		which is equivalent to 
		\begin{align*}
			\begin{aligned}
				&\ \langle\vartheta, w-u\rangle + J(w)-J(u)
				\leq \langle \vartheta - \eta, w - y\rangle + V(r, u)-V(y, u)-V(r,w)-V(w, y).
			\end{aligned}
		\end{align*}		
		Applying Schwartz inequality and Young's inequality to the above inequality, and using the well-known result that
		\begin{align}
			\label{eqnVvsNorm}
			\frac{\mu}{2}\|z- u\|^2\leq V(u, z), \forall u, z\in Z,
		\end{align}
		we obtain
		\begin{align}
			\label{tmp1}
			\begin{aligned}
				&\ \langle\vartheta, w-u\rangle + J(w)-J(u)\\
				\leq &\ \|\vartheta - \eta\|_* \|w - y\| + V(r, u)-V(y, u) -V(r, w)-\frac{\mu}2\|w-y\|^2
				\\
				\leq &\ \frac{1}{2\mu}\|\vartheta - \eta\|_*^2 + \frac{\mu}2\|w - y\|^2 + V(r, u)-V(y, u) - V(r, w)-\frac{\mu}2\|w-y\|^2
				\\
				= &\ \frac{1}{2\mu}\|\vartheta - \eta\|_*^2 + V(r, u)-V(y, u) - V(r, w).
			\end{aligned}
		\end{align}		
		The result in \eqref{eqnPRecursion} then follows immediately from above relation, \eqref{eqnVvsNorm} and \eqref{eqnLM}.
		
		Moreover, observe that by setting $u=w$ and $u=y$ in \eqref{eqnMP1} and \eqref{tmp1} respectively, we have
		\begin{align*}
			\langle\vartheta, y- w\rangle + J(y) - J(w) & \leq V(r, w) - V(r, y) - V(y, w),
			\\
			\langle\vartheta, w-y\rangle + J(w)-J(y) 
			& \leq \frac{1}{2\mu}\|\vartheta - \eta\|_*^2 + V(r, y)  - V(r, w).
		\end{align*}		
		Adding the two inequalities above, and using \eqref{eqnVvsNorm} and \eqref{eqnLM}, we have
		\[0\leq\frac{1}{2\mu}\|\vartheta - \eta\|_*^2-V(y,w)\leq \frac{L^2}{2\mu}\|r - w\|^2 + \frac{M^2}{2\mu} - V(y, w)\leq  \frac{L^2}{\mu^2}V(r,w) + \frac{M^2}{2\mu} - V(y, w), \]
		thus \eqref{eqnVvr} holds.
	\end{proof}
\end{pro}

\vgap 

\begin{pro}
	\label{propSimplifiedQ}
	For any sequences $\{r_t\}_{t\geq 1}$ and $\{w_t\}_{t\geq 1}\subset Z$, if the sequences $\{w\ag\}$ and $\{w\md\}$ are generated by \eqref{eqnwmd} and \eqref{eqnwag}, then for all $u\in Z$, 
		\begin{align}
			\label{eqnSimplifiedQ}
			\begin{aligned}
				&\ Q(w\ag[1], u) - (1-\alpha_t)Q(w\ag, u)\\
				\leq &\alpha_t\langle\nabla G(w\md),w_{t+1}-u\rangle
				+\frac {L_G\alpha_t^2}{2}\|w_{t+1}-r_t\|^2
				 + \alpha_t\langle H(w\tn), w_{t+1}-u\rangle +  \alpha_tJ(w_{t+1})  -\alpha_tJ(u).
			\end{aligned}
		\end{align}
	\begin{proof}
		Observing from (\ref{eqnwmd}) and (\ref{eqnwag}) that $w\ag[1]-w\md=\alpha_t(w_{t+1}-r_t)$. This observation together with the convexity of $G(\cdot)$, then imply that for all $u\in Z$,
		\begin{align*}
			 G(w\ag[1])
			\leq\ &  G(w\md)+\langle\nabla G(w\md), w\ag[1]-w\md\rangle + \frac {L_G}2\|w\ag[1]-w\md\|^2
			\\
			=\ & (1-\alpha_t)\left[G(w\md)+\langle\nabla G(w\md),w\ag-w\md\rangle\right]\\
				&+\alpha_t\left[G(w\md)+\langle\nabla G(w\md),u-w\md\rangle\right]+\alpha_t\langle\nabla G(w\md),w_{t+1}-u\rangle
				+\frac {L_G\alpha_t^2}{2}\|w_{t+1}-r_t\|^2
			\\
			\leq\ &(1-\alpha_t)G(w\ag)+\alpha_tG(u)+\alpha_t\langle\nabla G(w\md),w_{t+1}-u\rangle
			+\frac {L_G\alpha_t^2}{2}\|w_{t+1}-r_t\|^2.
		\end{align*}
		Applying \eqref{eqnwag} and \eqref{eqnQ} to the above inequality, and using the fact that $H(\cdot)$ is monotone, we have
		\begin{align*}
			\begin{aligned}
				&\ Q(w\ag[1], u) - (1-\alpha_t) Q(w\ag, u)\\
				=\ &\ G(w\ag[1])-(1-\alpha_t)G(w\ag)-\alpha_tG(u)
				+ \langle H(u), w\ag[1]-u\rangle -(1-\alpha_t)\langle H(u), w\ag-u\rangle \\
				& +  J(w\ag[1]) - (1-\alpha_t)J(w\ag) -\alpha_tJ(u)\\
				\leq\ &\ G(w\ag[1])-(1-\alpha_t)G(w\ag)-\alpha_tG(u)
				+ \alpha_t\langle H(u), w_{t+1}-u\rangle +  \alpha_tJ(w_{t+1})  -\alpha_tJ(u)
				\\
				\leq &\ \alpha_t\langle\nabla G(w\md),w_{t+1}-u\rangle
				+\frac {L_G\alpha_t^2}{2}\|w_{t+1}-r_t\|^2
				+ \alpha_t\langle H(w\tn), w_{t+1}-u\rangle +  \alpha_tJ(w_{t+1})  -\alpha_tJ(u).
			\end{aligned}
		\end{align*}
	\end{proof}
\end{pro}

\vgap

The following lemma estimates a bound on $Q(w\ag[1], u)$, and will be used in the proof of both Theorems \ref{thmAMPRateB} and \ref{thmAMPRateUB}.

\begin{lem}
	\label{lemQBoundGeneral}
	Suppose that the parameters $\{\alpha_t\}$ in Algorithm \ref{algAMP} satisfies $\alpha_1=1$ and $0\leq\alpha_t<1$ for all $t>1$, and let the sequence $\{\Gamma_t\}$ be defined in \eqref{eqnGamma}.
	Then the iterates $\{r_t\}, \{w_t\}$ and $ \{w\ag\}$ of Algorithm \ref{algAMP} satisfy
	\begin{align}
		\label{eqnQBoundGeneral}
		& \frac 1{\Gamma_t}Q(w\ag[1], u) \leq   
		{\cal B}_t(u, r_{[t]})
		-\sum_{i=1}^t\frac{\alpha_i}{2\Gamma_i\gamma_i}\left({\mu}-{L_G\alpha_i\gamma_i}- \frac{L_H^2\gamma_i^2}{\mu}\right)\|r_i-w_{i+1}\|^2, \ \forall u\in Z,
		\\
		\label{eqnB}
		\text{where }&{\cal B}_t(u, r_{[t]}):=\sum_{i=1}^t\frac{\alpha_i}{\Gamma_i\gamma_i}(V(r_i, u) - V(r_{i+1},u)).
	\end{align}	
	\begin{proof}
		First, it follows from Proposition \ref{propPRecursion} applied to iterations (\ref{eqnProxwmd}) and (\ref{eqnProxrmd}) with $r=r_t, w=w_{t+1}, y=r_{t+1}, \vartheta = \gamma_tH(r_t) + \gamma_t\nabla G(w\md), \eta=\gamma_tH(w\tn) + \gamma_t\nabla G(w\md), J=\gamma_tJ$, $L=L_H\gamma_t$ and $M=0$ that for any $u\in Z$,
		\begin{equation*}
			\begin{aligned}
				& \gamma_t\langle H(w\tn) + \nabla G(w\md), w_{t+1}-u\rangle + \gamma_tJ(w_{t+1})-\gamma_tJ(u)\\
				\leq &\ V(r_t, u)-V(r_{t+1}, u)-\left(\frac{\mu}2- \frac{L_H^2\gamma_t^2}{2\mu}\right)\|r_t-w_{t+1}\|^2.
			\end{aligned}
		\end{equation*}
		Now applying the above inequality to \eqref{eqnSimplifiedQ}, we have
		\begin{align}
			\begin{aligned}
				\ Q(w\ag[1], u) - (1-\alpha_t)Q(w\ag, u)
				\leq &\ \frac{\alpha_t}{\gamma_t}\left[V(r_t, u)-V(r_{t+1}, u)\right]-\frac{\alpha_t}{2\gamma_t}\left({\mu}-{L_G\alpha_t\gamma_t}- \frac{L_H^2\gamma_t^2}{\mu}\right)\|r_t-w_{t+1}\|^2.
			\end{aligned}
		\end{align}
		Dividing both sides of the above inequality by $\Gamma_t$, we have
		\begin{align*}
			\begin{aligned}
				&\ \frac 1{\Gamma_t}Q(w\ag[1], u) - \frac{1-\alpha_t}{\Gamma_t}Q(w\ag, u)
				\\
				\leq &\ \frac{\alpha_t}{\Gamma_t\gamma_t}\left[V(r_t, u)-V(r_{t+1}, u)\right]-\frac{\alpha_t}{2\Gamma_t\gamma_t}\left({\mu}-{L_G\alpha_t\gamma_t}- \frac{L_H^2\gamma_t^2}{\mu}\right)\|r_t-w_{t+1}\|^2.
			\end{aligned}
		\end{align*}		
		Using the facts that $\alpha_1=1$, and that $\ds\frac{1-\alpha_t}{\Gamma_t} = \frac 1{\Gamma_{t-1}},\ t>1,$ due to \eqref{eqnGamma}, we can apply the above inequality recursively to obtain \eqref{eqnQBoundGeneral}.
	\end{proof}
\end{lem}

\vgap

We are now ready to prove Theorem \ref{thmAMPRateB}, which provides an estimate of the gap function of deterministic AMP when $Z$ is bounded. This result follows immediately from Lemma \ref{lemQBoundGeneral}.

{\it Proof of Theorem \ref{thmAMPRateB}.
		In view of \eqref{eqnCondAlphaGammaBD} and  \eqref{eqnQBoundGeneral}, to prove \eqref{eqngBound} it suffices to show that
		${\cal B}_t(u, r_{[t]})\leq \ds\frac{\alpha_t}{\Gamma_t\gamma_t}\Omega_Z^2$ for all $u\in Z$. In fact, since the sequence $\{r_i \}_{i=1}^{t+1}$ is in the bounded set $Z$, applying \eqref{eqnVBounded} and \eqref{eqnCondAlphaGammaBD} to \eqref{eqnB} we have
		\begin{align}
		\label{eqnBBD}
		\begin{aligned}
			&\ {\cal B}_t(u, r_{[t]})
			=
			\frac{\alpha_1}{\Gamma_1\gamma_1}V(r_1, u) - \sum_{i=1}^{t-1}\left[\frac{\alpha_i}{\Gamma_i\gamma_i} - \frac{\alpha_{i+1}}{\Gamma_{i+1}\gamma_{i+1}} \right]V(r\tn[i], u) - \frac{\alpha_t}{\Gamma_t\gamma_t}V(r\tn, u)
			\\
			\leq &\ \frac{\alpha_1}{\Gamma_1\gamma_1}\Omega_Z^2 - \sum_{i=1}^{t-1}\left[\frac{\alpha_i}{\Gamma_i\gamma_i} - \frac{\alpha_{i+1}}{\Gamma_{i+1}\gamma_{i+1}} \right]\Omega_Z^2
			= \frac{\alpha_t}{\Gamma_t\gamma_t}\Omega_Z^2,\ \forall u\in Z
		\end{aligned}
		\end{align}		
		thus \eqref{eqngBound} holds.
	\endproof
}

\vgap 

In the remaining part of this subsection, we will focus on proving Theorem \ref{thmAMPRateUB}, which summarizes the convergence properties of deterministic AMP when $Z$ is unbounded.

{\it Proof of Theorem \ref{thmAMPRateUB}. 
	By the assumption that $V(r,z):=\|z-r\|^2/2$ for all $r,z\in Z$, and applying the last relation of \eqref{eqnCondAlphaGammaUB} to \eqref{eqnB}, we obtain
	\[{\cal B}_t(u, r_{[t]}) =  \frac{\alpha_t}{2\Gamma_t\gamma_t}\|r_1 - u\|^2
			 - \frac{\alpha_t}{2\Gamma_t\gamma_t}\|r\tn - u\|^2. \]
	Applying this and the second relation of \eqref{eqnCondAlphaGammaUB} to \eqref{eqnQBoundGeneral} and noting that $\mu=1$, we have
	\begin{align}
		\label{eqnQEucl}
		Q(w\ag[1], u) \leq &\   \frac{\alpha_t}{2\gamma_t}\|r_1 - u\|^2
		- \frac{\alpha_t}{2\gamma_t}\|r\tn - u\|^2
		-\frac{\alpha_t}{2\gamma_t}\sum_{i=1}^t\left({1}-c^2\right)\|r_i-w_{i+1}\|^2.
	\end{align}
	The first two terms in \eqref{eqnQEucl} can be rewritten as
	\begin{align}
		\label{eqnr2wag}
		\begin{aligned}
		&\ \frac 12\|r_1 - u\|^2 - \frac 12\|r\tn - u\|^2  = \frac 12 \|r_1\|^2 - \frac 12\|r\tn\|^2 - \langle r_1 - r\tn, u\rangle
		\\
		= &\ \frac 12\|r_1 - w\ag[1]\|^2 - \frac 12\|r\tn - w\ag[1]\|^2 + \langle r_1 - r\tn, w\ag[1] - u\rangle.
		\end{aligned}
	\end{align}
	Then, the combination of \eqref{eqnQEucl} and \eqref{eqnr2wag} yields
	\begin{align}
		\label{eqneps}
\begin{aligned}
		&\ Q(w\ag[1], u) - \frac{\alpha_t}{\gamma_t}\langle r_1 - r\tn, w\ag[1] - u\rangle \\
		\leq &\ \frac{\alpha_t}{2\gamma_t}\|r_1 - w\ag[1]\|^2 - \frac{\alpha_t}{2\gamma_t}\|r\tn - w\ag[1]\|^2 - \frac{\alpha_t}{2\gamma_t}(1-c^2)\sum_{i=1}^t\|r_i-w_{i+1}\|^2=:\varepsilon\tn.
\end{aligned}
	\end{align}
	Therefore, if we set $\ds v\tn : = \frac{\alpha_t}{\gamma_t}(r_1 - r\tn),$	then 
	$Q(w\ag[1], u) - \langle v\tn, w\ag[1] - u\rangle \leq \varepsilon\tn$ for all $u\in Z$. It should be noted that $\varepsilon\tn\geq 0$ holds trivially by letting $u=w\ag[1]$ in \eqref{eqneps}. Hence we have $\tilde g(w\ag[1], v\tn)\leq \varepsilon\tn$ and it suffices to estimate the bound of $\|v\tn\|$ and $\varepsilon\tn$. 

	If there exists a strong solution $u^*$ of $VI(Z;G,H,J)$, then by  \eqref{eqnF}, \eqref{eqnSVI}, \eqref{eqnQ} and the convexity of $G$ and $J$, we have $Q(w\ag[1], u^*)\geq \langle \nabla F(u^*),w\ag[1]-u^*\rangle\ge 0$. This observation together with  \eqref{eqnQEucl} imply that 
		\begin{align*}
			\|r_1 - u^*\|^2
			- \|r\tn - u^*\|^2
			-\sum_{i=1}^t\left({1}-c^2\right)\|r_i-w_{i+1}\|^2 \geq 0.
		\end{align*}
		By the above inequality and the definition of $D$ in \eqref{eqndtheta}, we have the following two inequalities:
		\begin{align}
			\label{eqnrtnBound}
			\|r\tn - u^*\|\leq & D, \\
			\label{eqnrtwtnBound}
			\sum_{i=1}^t\|r_i-w_{i+1}\|^2 \leq & \frac{D^2}{1 - c^2}.
		\end{align}
	By \eqref{eqnrtnBound} and the definition of $v\tn$, we have
	\begin{align*}
		\|v\tn\| \leq \frac{\alpha_t}{\gamma_t}\left(\|r_1 - u^*\| + \|r\tn - u^*\|\right) \leq \frac{2\alpha_t}{\gamma_t}D,
	\end{align*}
	hence the first relation in \eqref{eqnvepsThm} holds.
	
	To finish the proof, it now suffices to estimate a bound for $\varepsilon_t$. Firstly we explore the definition of the aggregate point $w\ag[1]$. By \eqref{eqnwag} and \eqref{eqnGamma}, we have
	\[\frac 1{\Gamma_t}w\ag[1] = \frac 1{\Gamma_{t-1}} w\ag + \frac{\alpha_t}{\Gamma_t}w\tn, \ \forall t\ge 1.\]
	Using the assumption that $w_1^{ag} = w_1$, we obtain
	\begin{align}
		\label{eqnwagReform}
		w\ag[1] = \Gamma_t\sum_{i=1}^t\frac{\alpha_i}{\Gamma_i}w_{i+1},
	\end{align}
	where by \eqref{eqnGamma} we have 
	\begin{align}
		\label{eqnGammaSpan}
		\Gamma_t\sum_{i=1}^t\frac{\alpha_i}{\Gamma_i} = 1.	
	\end{align}  
	Therefore, $w\ag[1]$ is a convex combination of iterates $w_2, \ldots, w\tn$. Using \eqref{eqndtheta}, \eqref{eqneps}, \eqref{eqnrtnBound} and \eqref{eqnrtwtnBound}, we conclude that
	\begin{align*}
		&\ \varepsilon\tn \le \frac{\alpha_t}{2\gamma_t}\|r_1-w\ag[1]\|^2\le \frac{\alpha_t\Gamma_t}{2\gamma_t}\sumt\frac{\alpha_i}{\gamma_i}\|r_1 - w\tn[i]\|^2 
		\\
		\le&\  \frac{3\alpha_t\Gamma_t}{2\gamma_t}\sumt\frac{\alpha_i}{\gamma_i}(\|r_1 - u^*\|^2 + \|r_i-u^*\|^2 + \|r_i- w\tn[i]\|^2) \le \frac{3\alpha_t}{2\gamma_t}\left(2D^2 + \Gamma_t\max_{i=1,\ldots,t}\frac{\alpha_i}{\gamma_i}\sumt\|r_i-w\tn[i]\|^2 \right)
		\\
		\le&\ \frac{3\alpha_t(1+\theta_t)D^2}{\gamma_t}.
	\end{align*}
	\endproof
}

\subsection{Convergence analysis for stochastic AMP}
\label{secProofS}
In this section, we prove the convergence results of the stochastic AMP method presented in Section \ref{secAMPS}, namely, Theorems \ref{thmAMPRateBS} and \ref{thmAMPRateUBS}.

Throughout this section, we will use the following notations to describe the inexactness of the first order information from ${\cal SO}_H$ and ${\cal SO}_G$. At the $t$-th iteration, letting $\cHrt$, $\cHwt$ and $\cGt$ be the outputs of the stochastic oracles, we denote 
		\begin{align}
			\label{eqnDelta}
			\Delta_H^{2t-1} := \cHrt - H(r_t),\ \Delta_H^{2t} := \cHwt - H(w\tn)\text{ and } \Delta_G^{t} := \cGt - \nabla G(w\md).
		\end{align}

To start with, we present a technical result to obtain a bound on $Q(w\ag[1], u)$ for all $u\in Z$. The following lemma is analogous to Lemma \ref{lemQBoundGeneral} for deterministic AMP, and will be applied in the proof of Theorems \ref{thmAMPRateBS} and \ref{thmAMPRateUBS}.
\begin{lem}
	\label{lemQBoundGeneralS}
	Suppose that the parameters $\{\alpha_t\}$ in Algorithm \ref{algAMP} satisfies $\alpha_1=1$ and $0\leq\alpha_t<1$ for all $t>1$, and let the sequence $\{\Gamma_t\}$ be defined in \eqref{eqnGamma}.
	Then the iterates $\{r_t\}, \{w_t\}$ and $\{w\ag\}$ generated by Algorithm \ref{algAMPS} satisfy
	\begin{align}
		\label{eqnQBoundGeneralS}
		& \frac 1{\Gamma_t}Q(w\ag[1], u) \leq   
		{\cal B}_t(u, r_{[t]})
		-\sum_{i=1}^t\frac{\alpha_i}{2\Gamma_i\gamma_i}\left(q{\mu}-{L_G\alpha_i\gamma_i}- \frac{3L_H^2\gamma_i^2}{\mu}\right)\|r_i-w_{i+1}\|^2 + \sum_{i=1}^t\Lambda_i(u), \ \forall u\in Z,
	\end{align}		
		where ${\cal B}_t(u, r_{[t]})$ is defined in \eqref{eqnB}, and
		\begin{equation}
			\label{eqnLambda}
			\Lambda_i(u):= \frac{3\alpha_i\gamma_i}{2\mu\Gamma_i}\left(\|\Delta_H^{2i}\|_*^2 + \|\Delta_H^{2i-1}\|_*^2 \right)  
			- \frac{(1-q)\mu\alpha_i}{2\Gamma_i\gamma_i}\|r_i - w_{i+1}\|^2
			- \frac{\alpha_i}{\Gamma_i}\langle \Delta_H^{2i}+\Delta_G^i, w_{i+1} - u\rangle.
		\end{equation}

	\begin{proof}
		Observe from \eqref{eqnDelta} that
		\begin{align}
			\label{eqnLMStoc}
			\begin{aligned}
				&\ \| \cHwt - \cHrt\|_*^2
				\leq \left(\| H(w\tn) - H(r_t)\|_* + \|\Delta_H^{2t}\|_* + \|\Delta_H^{2t-1}\|_*\right)^2
				\\
				\leq &\ 3 \left(\| H(w\tn) - H(r_t)\|_*^2 + \|\Delta_H^{2t}\|_*^2 + \|\Delta_H^{2t-1}\|_*^2\right)
				\\
				\leq &\ 3 \left(L_H^2\|w\tn - r_t\|^2 + \|\Delta_H^{2t}\|_*^2 + \|\Delta_H^{2t-1}\|_*^2\right).
			\end{aligned}
		\end{align}		
		Applying Proposition \ref{propPRecursion} to \eqref{eqnProxwmdS} and \eqref{eqnProxrmdS} with $r=r_t, w=w_{t+1}, y=r_{t+1}, \vartheta = \gamma_t\cHrt + \gamma_t\cGt,$ $\eta=\gamma_t\cHwt + \gamma_t\cGt, J=\gamma_tJ$, $L^2=3L_H^2\gamma_t^2$ and $M^2 = 3\gamma_t^2(\|\Delta_H^{2t}\|_*^2 + \|\Delta_H^{2t-1}\|_*^2)$, and using \eqref{eqnLMStoc} we get that for any $u\in Z$,
		\begin{align*}
			\begin{aligned}
				&\ \gamma_t\langle \cHwt + \cGt, w\tn-u\rangle  + \gamma_t J(w)-\gamma_t J(u)
				\\
				\leq &\ V(r_t, u)-V(r\tn, u) - \left(\frac \mu 2 - \frac{3L_H^2\gamma_t^2}{2\mu}\right)\|r_t - w\tn\|^2 + \frac{3\gamma_t^2}{2\mu} (\|\Delta_H^{2t}\|_*^2 + \|\Delta_H^{2t-1}\|_*^2).
			\end{aligned}
		\end{align*}
		Applying \eqref{eqnDelta} and the above inequality to \eqref{eqnSimplifiedQ}, we have
		\begin{align*}
			\begin{aligned}
				&\ Q(w\ag[1], u) - (1-\alpha_t)Q(w\ag, u)
				\\
				\le &\ \alpha_t\langle \cHwt + \cGt, w_{t+1}-u\rangle 
				+  \alpha_tJ(w_{t+1})  -\alpha_tJ(u) +\frac {L_G\alpha_t^2}{2}\|w_{t+1}-r_t\|^2\\
				&\ - \alpha_t\langle \Delta_H^{2t} + \Delta_G^t, w\tn - u\rangle 
				\\			
				\leq &\ \frac{\alpha_t}{\gamma_t}\left(V(r_t, u) - V(r\tn, u) \right)  - \frac{\alpha_t}{2\gamma_t}\left(\mu - {L_G\alpha_t\gamma_t} -  \frac{3L_H^2\gamma_t^2}{\mu} \right)\|r_t - w\tn\|^2 + \frac{3\alpha_t\gamma_t}{2\mu}\left(\|\Delta_H^{2t}\|_*^2 + \|\Delta_H^{2t-1}\|_*^2 \right)
				\\
				&\ - \alpha_t\langle \Delta_H^{2t}+\Delta_G^t, w\tn - u\rangle.
			\end{aligned}
		\end{align*}
		Dividing the above inequality by $\Gamma_t$ and using the definition of $\Lambda_t(u)$ in \eqref{eqnLambda}, we obtain
		\begin{align*}
			\begin{aligned}
				&\ \frac{1}{\Gamma_t}Q(w\ag[1], u) - \frac{1-\alpha_t}{\Gamma_t}Q(w\ag, u)
				\\
				\leq &\ \frac{\alpha_t}{\Gamma_t\gamma_t}\left(V(r_t, u) - V(r\tn, u) \right)  - \frac{\alpha_t}{2\Gamma_t\gamma_t}\left(q\mu - {L_G\alpha_t\gamma_t} -  \frac{3L_H^2\gamma_t^2}{\mu} \right)\|r_t - w\tn\|^2 + \Lambda_t(u).
			\end{aligned}
		\end{align*}		
		Noting the fact that $\alpha_1=1$ and $\ds\frac{1-\alpha_t}{\Gamma_t} = \frac 1{\Gamma_{t-1}}$, $t>1$, due to \eqref{eqnGamma}, applying the above inequality recursively and using the definition of  $\cB_t(\cdot,\cdot)$ in \eqref{eqnB}, we conclude \eqref{eqnQBoundGeneralS}.
		\end{proof}
	\end{lem}

\vgap 

We still need the following technical result to prove Theorem \ref{thmAMPRateBS}.

\begin{lem}
	\label{lemTech}
	Suppose that the sequences $\{\theta_t \}$ and $\{\gamma_t \}$ are positive sequences. For any $w_1\in Z$ and any sequence $\{\Delta^t \}\subset\cE$, if we define $w^v_1=w_1$ and 
	\begin{equation}
		\label{eqnProxv}
		w_{i+1}^v = \argmin_{u\in Z}-\gamma_i\langle\Delta^i, u\rangle + V(w_i^v, u), \ \forall i>1,
	\end{equation}
	then
	\begin{equation}
		\label{eqnTech}
		\sum_{i=1}^t\theta_i\langle-\Delta^i, w_i^v - u\rangle \leq \sumt\frac{\theta_i}{\gamma_i}(V(w_i^v, u)  - V(w_{i+1}^v, u)) + \sum_{i=1}^t\frac{\theta_i\gamma_i}{2\mu}\|\Delta_i\|_*^2, \ \forall u\in Z.
	\end{equation}
	\begin{proof}
		Applying Proposition \ref{propProxMap} with $r=w_i^v$, $w=w_{i+1}^v$, $\zeta=-\gamma_i\Delta^i$ and $J=0$, we have
		\begin{align*}
			-\gamma_i\langle\Delta^i, w_{i+1}^v - u\rangle \leq V(w_i^v, u) -V(w_i^v, w_{i+1}^v)  - V(w_{i+1}^v, u), \ \forall u\in Z.
		\end{align*}
		On the other hand, by Schwartz inequality, Young's inequality and \eqref{eqnVvsNorm} we have
		\begin{align*}
			-\gamma_i\langle\Delta^i, w_i^v - w_{i+1}^v\rangle \leq \gamma_i\|\Delta^i\|_*\|\|w_i^v - w_{i+1}^v\| \leq\frac{\gamma_i^2}{2\mu}\|\Delta_i\|_*^2 +  \frac{\mu}{2}\|w_i^v - w_{i+1}^v\|^2 \leq \frac{\gamma_i^2}{2\mu}\|\Delta_i\|_*^2 + V(w_i^v, w_{i+1}^v).
		\end{align*}
		Adding the two inequalities above and multiplying them by $\theta_i/\gamma_i$, we obtain
		\begin{align*}
			-\theta_i\langle\Delta^i, w_i^v-u\rangle \leq \frac{\theta_i\gamma_i}{2\mu}\|\Delta_i\|_*^2 + \frac{\theta_i}{\gamma_i}(V(w_i^v, u)  - V(w_{i+1}^v, u)).
		\end{align*}		
		Summing from $i=1$ to $t$, we conclude \eqref{eqnTech}.
	\end{proof}
\end{lem}

\vgap 

We are now ready to prove Theorem \ref{thmAMPRateBS}.

{\it Proof of Theorem \ref{thmAMPRateBS}.
	Firstly, applying \eqref{eqnCondAlphaGammaIncS} and \eqref{eqnBBD} to \eqref{eqnQBoundGeneralS} in Lemma \ref{lemQBoundGeneralS}, we have
	\begin{align}
		\label{eqnQLambda}
		\frac{1}{\Gamma_t}Q(w\ag[1], u) \leq \frac{\alpha_t}{\Gamma_t\gamma_t}\Omega_Z^2 + \sum_{i=1}^{t}\Lambda_i(u), \ \forall u\in Z.
	\end{align}
	
	Letting $w_1^v = w_1$, defining $w_{i+1}^v$ as in \eqref{eqnProxv} with $\Delta^i = \Delta^{2i}_H + \Delta^i_G$ for all $i>1$, it follows from \eqref{eqnB} and Lemma \ref{lemTech} with  $\theta_i=\alpha_i/\Gamma_i$ that
	\begin{align}
		\label{tmp3}
		-\sumt\frac{\alpha_i}{\Gamma_i}\langle \Delta^{2i}_H + \Delta^i_G, w^v_i-u\rangle \le {\cal B}_t(u, w^v_{[t]}) +  \sum_{i=1}^t\frac{\alpha_i\gamma_i}{2\mu\Gamma_i}\|\Delta_H^{2i} + \Delta_G^i\|_*^2,\ \forall u\in Z.
	\end{align}
	Noting that by \eqref{eqnLambda}
	\begin{align*}
		\begin{aligned}
		\sum_{i=1}^{t}\Lambda_i(u) = &\  - \sum_{i=1}^{t}\frac{\alpha_i}{\Gamma_i}\langle \Delta_H^{2i} + \Delta_G^i, w_{i}^v - u\rangle	+ \sum_{i=1}^t\frac{\alpha_i}{\Gamma_i}\left[-\frac{(1-q)\mu}{2\gamma_i}\|r_i - w_{i+1}\|^2 - \langle \Delta_G^i, w_{i+1} - r_i\rangle \right] 
		\\
		&\ + \sum_{i=1}^{t}\frac{3\alpha_i\gamma_i}{2\mu\Gamma_i}\left(\|\Delta_H^{2i}\|_*^2 + \|\Delta_H^{2i-1}\|_*^2 \right)	
		- \sum_{i=1}^{t}\frac{\alpha_i}{\Gamma_i}\langle \Delta_G^i, r_i - w_i^v\rangle
		- \sum_{i=1}^{t}\frac{\alpha_i}{\Gamma_i}\langle \Delta_H^{2i}, w_{i+1} - w_{i}^v\rangle,
		\end{aligned}
	\end{align*}
	applying \eqref{tmp3} and the Young's inequality to above equation, we conclude that
	\begin{align}
		\label{eqnLambdaSimplifed}
		\begin{aligned}
		&\ \sum_{i=1}^{t}\Lambda_i(u) \leq {\cal B}_t(u, w^v_{[t]}) + U_t,
		\end{aligned}
	\end{align}
	where
	\begin{align}
		\label{eqnU}	
		\begin{aligned}
		 U_t := &\  \sum_{i=1}^t\frac{\alpha_i\gamma_i}{2\mu\Gamma_i}\|\Delta_H^{2i} + \Delta_G^i\|_*^2 + \sum_{i=1}^t\frac{\alpha_i\gamma_i}{2(1-q)\mu\Gamma_i}\|\Delta_G^i\|_*^2 + \sum_{i=1}^{t}\frac{3\alpha_i\gamma_i}{2\mu\Gamma_i}\left(\|\Delta_H^{2i}\|_*^2 + \|\Delta_H^{2i-1}\|_*^2 \right)	
		 \\
		&\ - \sum_{i=1}^{t}\frac{\alpha_i}{\Gamma_i}\langle \Delta_G^i, r_i - w_i^v\rangle
		- \sum_{i=1}^{t}\frac{\alpha_i}{\Gamma_i}\langle \Delta_H^{2i}, w_{i+1} - w_{i}^v\rangle.
		\end{aligned}
	\end{align}
	
	Applying \eqref{eqnBBD} and \eqref{eqnLambdaSimplifed} to \eqref{eqnQLambda}, we have
	\[\frac{1}{\Gamma_t}Q(w\ag[1], u) \leq \frac{2\alpha_t}{\gamma_t\Gamma_t}\Omega_Z^2 + U_t, \ \forall u\in Z, \]
	or equivalently,
	\begin{align}
		\label{eqngU}
		g(w\ag)\leq \frac{2\alpha_t}{\gamma_t}\Omega_Z^2 + \Gamma_tU_t.	
	\end{align}
	Now it suffices to bound $U_t$, in both expectation and probability.
	
	We prove part (a) first. By our assumptions on $\cSO_G$ and $\cSO_H$ and in view of \eqref{eqnProxwmdS}, \eqref{eqnProxrmdS} and \eqref{eqnProxv}, during the $i$-th iteration of Algorithm \ref{algAMPS}, the random noise $\Delta^{2i}_H$ is independent of $w_{i+1}$ and $w_i^v$, and $\Delta^{i}_G$ is independent of $r_i$ and $w_i^v$, hence $\E[\langle \Delta_G^i, r_i - w_i^v\rangle] = \E[\langle \Delta_H^{2i}, w_{i+1} - w_{i}^v\rangle] = 0 $. In addition, Assumption \ref{itmVB} implies that $\E[\|\Delta_G^i\|_*^2]\leq\sigma_{G}^2$, $\E[\|\Delta_H^{2i-1}\|_*^2]\leq\sigma_H^2$ and $\E[\|\Delta_H^{2i}\|_*^2]\leq\sigma_H^2$, where $\Delta_G^i$, $\Delta_H^{2i-1}$ and $\Delta_H^{2i}$ are independent. Therefore, taking expectation on \eqref{eqnU} we have
	\begin{align}
		\label{eqnEUt}
		\begin{aligned}
		\E[U_t] \leq&\ \E\left[ \sum_{i=1}^t\frac{\alpha_i\gamma_i}{\mu\Gamma_i}\left(\|\Delta_H^{2i}\|^2 + \|\Delta_G^i\|_*^2\right) + \sum_{i=1}^t\frac{\alpha_i\gamma_i}{2(1-q)\mu\Gamma_i}\|\Delta_G^i\|_*^2 + \sum_{i=1}^{t}\frac{3\alpha_i\gamma_i}{2\mu\Gamma_i}\left(\|\Delta_H^{2i}\|_*^2 + \|\Delta_H^{2i-1}\|_*^2 \right)\right]
		\\
		=&\ \sum_{i=1}^{t}\frac{\alpha_i\gamma_i}{\mu\Gamma_i}\left[ 4\sigma_H^2 + \left(1 + \frac{1}{2(1-q)} \right)\sigma_G^2\right].
		\end{aligned}
	\end{align}
	Taking expectation on both sides of \eqref{eqngU}, and using above estimation on $\E[U_t]$, we obtain \eqref{eqnQBoundBS}.
	
	Next we prove part (b). Observing that the sequence $\{\langle \Delta_G^i, r_i - w_i^v\rangle \}_{i\geq 1}$ is a martingale difference and hence satisfies the large-deviation theorem (see, e.g., Lemma 2 of \cite{lan2012validation}), therefore using Assumption \ref{itmLT} and the fact that
	\begin{align*}
		&\ \E[\exp\{\mu(\alpha_i\Gamma_i^{-1}\langle \Delta_G^i, r_i - w_i^v\rangle)^2/2(\sigma_G\alpha_i\Gamma_i^{-1}\Omega_Z)^2 \}] \leq \E[\exp\{\mu\|\Delta_G^i\|_*^2\|r_i-w_i^v\|^2/2\sigma_G^2\Omega_Z^2 \} ]
		\\
		\leq &\ \E[\exp\{\|\Delta_G^i\|_*^2 \}/\sigma_G^2 ] \leq \exp\{ 1\},
	\end{align*}
	we conclude from the large-deviation theorem that
	\begin{align}
	\label{tmpp1}
	Prob\left\{\sum_{i=1}^t\frac{\alpha_i}{\Gamma_i}\langle \Delta_G^i, r_i - w_i^v\rangle > \lambda\sigma_G\Omega_Z\sqrt{\frac{2}{\mu}\sum_{i=1}^t\left(\frac{\alpha_i}{\Gamma_i} \right)^2 }  \right\}  \leq \exp\{-\lambda^2/3\}.
	\end{align}
	By similar argument we also have
	\begin{align}
	\label{tmpp2}
	Prob\left\{\sum_{i=1}^t\frac{\alpha_i}{\Gamma_i}\langle \Delta_H^{2i}, w_{i+1} - w_i^v\rangle > \lambda\sigma_H\Omega_Z\sqrt{\frac{2}{\mu}\sum_{i=1}^t\left(\frac{\alpha_i}{\Gamma_i} \right)^2 }  \right\}  \leq \exp\{-\lambda^2/3\}.
	\end{align}
	
	In addition,  letting $S_i = \alpha_i\gamma_i/(\mu\Gamma_i)$ and $S = \sum_{i=1}^tS_i$, by Assumption \ref{itmLT} and the convexity of exponential functions, we have
	\[
			\begin{aligned}
				\E\left[\exp\left\{\frac 1S\sum_{i=1}^{t}S_i{\|\Delta_{G}^i\|_*^2}/{\sigma_{G}^2} \right\}\right] \leq \E\left[\frac 1S\sum_{i=1}^{t}S_i\exp\left\{\|\Delta_{G}^i\|_*^2/\sigma_{G}^2 \right\}\right]\leq\exp\{1\},
			\end{aligned}
	\]	
	therefore by Markov's inequality we have
	\begin{align}
		\label{tmpp3}
		Prob\left\{\left(1 + \frac{1}{2(1-q)} \right)\sum_{i=1}^t \frac{\alpha_i\gamma_i}{\mu\Gamma_i}\|\Delta^i_G\|_*^2 > (1+\lambda)\sigma_G^2\left(1 + \frac{1}{2(1-q)} \right)\sum_{i=1}^t \frac{\alpha_i\gamma_i}{\mu\Gamma_i} \right\} \leq \exp\{-\lambda\}.
	\end{align}
	Using similar arguments, we also have
	\begin{align}
		\label{tmpp4}
		& Prob\left\{\sum_{i=1}^t \frac{3\alpha_i\gamma_i}{2\mu\Gamma_i}\|\Delta^{2i-1}_H\|_*^2 > (1+\lambda)\frac{3\sigma_H^2}{2}\sum_{i=1}^t \frac{\alpha_i\gamma_i}{\mu\Gamma_i} \right\} \leq \exp\{-\lambda\},
		\\
		\label{tmpp5}
		& Prob\left\{\sum_{i=1}^t \frac{5\alpha_i\gamma_i}{2\mu\Gamma_i}\|\Delta^{2i}_H\|_*^2 > (1+\lambda)\frac{5\sigma_H^2}{2}\sum_{i=1}^t \frac{\alpha_i\gamma_i}{\mu\Gamma_i} \right\} \leq \exp\{-\lambda\}.
	\end{align}		
	Using the fact that $\|\Delta_H^{2i} + \Delta_G^{2i-1}\|_{*}^2\le 2\|\Delta_H^{2i}\|_{*}^2 + \|\Delta_G^{2i-1}\|_{*}^2$, we conclude from \eqref{eqngU}--\eqref{tmpp5} that \eqref{eqnProb} holds.
\endproof}

\vgap 

In the remaining part of this subsection, we will focus on proving Theorem \ref{thmAMPRateUBS}, which describes the rate of convergence of Algorithm \ref{algAMPS} for solving $SVI(Z;G,H,J)$ when $Z$ is unbounded.

{\it Proof the Theorem \ref{thmAMPRateUBS}.
Let $U_t$ be defined in \eqref{eqnU}. Firstly, applying \eqref{eqnCondAlphaGammaUBS} and \eqref{eqnLambdaSimplifed} to \eqref{eqnQBoundGeneralS} in Lemma \ref{lemQBoundGeneralS}, we have
	\begin{align}
		\label{eqnQEuclS}
		& \frac 1{\Gamma_t}Q(w\ag[1], u) \leq   
		{\cal B}_t(u, r_{[t]})
		-\frac{\alpha_t}{2\Gamma_t\gamma_t}\sum_{i=1}^t\left(q-c^2\right)\|r_i-w_{i+1}\|^2 + {\cal B}_t(u, w^v_{[t]}) + U_t, \ \forall u\in Z.
	\end{align}		
In addition, applying \eqref{eqnCondAlphaGammaUBS} to the definition of $\cB_t(\cdot,\cdot)$ in \eqref{eqnB}, we obtain
	\begin{align}
	\label{eqnBr}
	 {\cal B}_t(u, r_{[t]}) 
	= &\ \frac{\alpha_t}{2\Gamma_t\gamma_t}(\|r_1 - u\|^2 - \|r\tn-u\|^2)
	\\
	\label{eqnBrw}
	= &\ \frac{\alpha_t}{2\Gamma_t\gamma_t}(\|r_1 - w\ag[1]\|^2 - \|r\tn-w\ag[1]\|^2 + 2\langle r_1 - r\tn, w\ag[1] - u\rangle).
	\end{align}
By using similar argument and the fact that $w^v_1=w_1=r_1$, we have
\begin{align}
	\label{eqnBv}
	 {\cal B}_t(u, w^v_{[t]}) 
	=&\ \frac{\alpha_t}{2\Gamma_t\gamma_t}(\|r_1 - u\|^2 - \|w^v_{t+1}-u\|^2)
	\\	
	\label{eqnBvw}
	=&\ \frac{\alpha_t}{2\Gamma_t\gamma_t}(\|r_1 - w\ag[1]\|^2 - \|w^v_{t+1}-w\ag[1]\|^2 + 2\langle r_1 - w_{t+1}^v, w\ag[1] - u\rangle).
\end{align}
We then conclude from \eqref{eqnQEuclS}, \eqref{eqnBrw} and \eqref{eqnBvw} that
\begin{align}
	\label{tmp}
	Q(w\ag[1], u) - \langle v\tn, w\ag[1] - u\rangle \leq \varepsilon\tn,\ \forall u\in Z,
\end{align}
where
\begin{align}
	\label{eqnvS}
	v\tn : = &\frac{\alpha_t}{\gamma_t}(2r_1 - r\tn - w_{t+1}^v),\text{ and }
	\\
	\label{eqnepsS}
	\varepsilon\tn : = &\frac{\alpha_t}{2\gamma_t}\left(2\|r_1 - w\ag[1]\|^2 - \|r\tn - w\ag[1]\|^2 - \|w_{t+1}^v - w\ag[1]\|^2 -\sum_{i=1}^t\left(q-c^2\right)\|r_i-w_{i+1}\|^2 \right) +\Gamma_tU_t.
\end{align}
It is easy to see that the residual $\varepsilon\tn$ is positive by setting $u=w\ag[1]$ in \eqref{tmp}. Hence $\tilde g(w\ag[1], v\tn) \leq\varepsilon\tn$. To finish the proof, it suffices to estimate the bounds for $\E[\|v\tn\|]$ and $\E[\varepsilon\tn]$.

Since $Q(w\ag[1], u^*)\geq 0$, letting $u=u^*$ in \eqref{eqnQEuclS}, we conclude from \eqref{eqnBr} and \eqref{eqnBv} that
\begin{align*}
	2\|r_1 - u^*\|^2 - \|r\tn - u^*\|^2 - \|w^v_{t+1} - u^*\|^2 - \sum_{i=1}^t\left(q-c^2\right)\|r_i-w_{i+1}\|^2 + \frac{2\Gamma_t\gamma_t}{\alpha_t}U_t \geq 0,
\end{align*}
and using the definition of $D$ in \eqref{eqndtheta}, we have
\begin{align}
	\label{eqnrwvBound}
	\|r\tn - u^*\|^2 + \|w^v_{t+1} - u^*\|^2 + \sum_{i=1}^t\left(q-c^2\right)\|r_i-w_{i+1}\|^2 \leq 2D^2 + \frac{2\Gamma_t\gamma_t}{\alpha_t}U_t.
\end{align}

In addition, applying \eqref{eqnCondAlphaGammaUBS} and the definition of $C_t$ in \eqref{eqnCtheta} to \eqref{eqnEUt}, we have
\begin{align}
	\label{eqnEUC}
	\E[U_t]\leq \sum_{i=1}^{t}\frac{\alpha_t\gamma_i^2}{\Gamma_t\gamma_t}\left[ 4\sigma_H^2 + \left(1 + \frac{1}{2(1-q)} \right)\sigma_G^2\right] = \frac{\alpha_t}{\Gamma_t\gamma_t}C_t^2.
\end{align}
Combining \eqref{eqnEUC} and \eqref{eqnrwvBound}, we have 
\begin{align}
	\label{eqnDC}
	&\ \E[\|r\tn - u^*\|^2] + \E[\|w^v_{t+1} - u^*\|^2] + \sum_{i=1}^t\left(q-c^2\right)\E[\|r_i-w_{i+1}\|^2]
	\leq 2D^2  + 2C_t^2.
\end{align}
We are now ready to prove \eqref{eqnEv}. Observing from the definition of $v\tn$ \eqref{eqnvS} and the definition of $D$ in \eqref{eqndtheta} that $\|v\tn\| \le \frac{\alpha_t}{\gamma_t}(2D + \|w\tn^v - u^*\| + \|r\tn - u^*\|)$, applying Jensen's inequality and \eqref{eqnDC}, we obtain
\begin{align*}
	&\ \E[\|v\tn\|] \le \frac{\alpha_t}{\gamma_t}(2D + \sqrt{\E[(\|r\tn - u^*\| + \|w\tn^v - u^*\|)^2]}) 
	\\
	\le&\ \frac{\alpha_t}{\gamma_t}(2D + \sqrt{2\E[\|r\tn - u^*\|^2 + \|w\tn^v - u^*\|^2] }) \le \frac{\alpha_t}{\gamma_t}(2D + 2\sqrt{D^2 + C_t^2}).
\end{align*}

Our remaining goal is to prove \eqref{eqnEeps}. By applying Proposition \ref{propPRecursion} to \eqref{eqnProxwmdS} and \eqref{eqnProxrmdS} with $r=r_t, w=w_{t+1}, y=r_{t+1}, \vartheta = \gamma_t\cHrt + \gamma_t\cGt, \eta=\gamma_t\cHwt + \gamma_t\cGt, J=\gamma_tJ$, $L=3L_H^2\gamma_t^2$ and $M^2 = 3\gamma_t^2(\|\Delta_H^{2t}\|_*^2 + \|\Delta_H^{2t-1}\|_*^2)$ and using \eqref{eqnLMStoc} and \eqref{eqnVvr}, we have
\begin{align*}
	\frac 12\|r\tn - w\tn\|^2 
	\leq&\ \frac{3L_H^2\gamma_t^2}{2}\|r_t - w\tn\|^2 + \frac{3\gamma_t^2}{2}(\|\Delta_H^{2t}\|_*^2 + \|\Delta_H^{2t-1}\|_*^2) 
	\\
	\leq&\ \frac{c^2}{2}\|r_t - w\tn\|^2 + \frac{3\gamma_t^2}{2}(\|\Delta_H^{2t}\|_*^2 + \|\Delta_H^{2t-1}\|_*^2),
\end{align*}
where the last inequality is from \eqref{eqnCondAlphaGammaUBS}. Now using \eqref{eqnwagReform}, \eqref{eqnGammaSpan}, \eqref{eqnepsS}, the inequality above, and applying Jensen's inequality, we have
	\begin{align}
		\begin{aligned}
		&\ \varepsilon\tn - \Gamma_tU_t \leq \frac{\alpha_t}{\gamma_t}\|r_1 - w\ag[1]\|^2 = \frac{\alpha_t}{\gamma_t}\left\|r_1 - u^* +  \sumt\frac{\alpha_i}{\Gamma_i}(u^*-r\tn[i]) + \sumt\frac{\alpha_i}{\Gamma_i}(r\tn[i] - w\tn[i])\right\|
		\\
		\leq &\ \frac{3\alpha_t}{\gamma_t}\left[D^2 + \Gamma_t\sum_{i=1}^t\frac{\alpha_i}{\Gamma_i}\left(\|r_{i+1} - u^*\|^2 + \|w_{i+1} - r_{i+1}\|^2\right)\right]
		\\
		\leq &\ \frac{3\alpha_t}{\gamma_t}\left[D^2 + \Gamma_t\sum_{i=1}^t\frac{\alpha_i}{\Gamma_i}\left(\|r_{i+1} - u^*\|^2 + c^2\|w_{i+1} - r_{i}\|^2 + 3\gamma_i^2(\|\Delta_H^{2i}\|_*^2 + \|\Delta_H^{2i-1}\|_*^2)\|\right)\right].
	\end{aligned}
	\label{tmp2}
	\end{align}	
	Noting that by \eqref{eqnCtheta} and \eqref{eqnrwvBound}
	\begin{align*}
		&\ \Gamma_t\sumt\aGi(\|r_{i+1} - u^*\|^2 + c^2\|w_{i+1} - r_{i}\|^2) 
		\le\ \Gamma_t\sumt\frac{\alpha_i\theta}{\Gamma_i}(\|r_{i+1} - u^*\|^2 + (q-c^2)\|w_{i+1} - r_{i}\|^2)
		\\
		\le &\ \Gamma_t\sumt\frac{\alpha_i\theta}{\Gamma_i}(2D^2 + \frac{2\Gamma_i\gamma_i}{\alpha_i}U_i)  = 2\theta D^2 + 2\theta\Gamma_t\sumt{\gamma_i}U_i,
	\end{align*}
	and that by \eqref{eqnCondAlphaGammaUBS}
	\begin{align*}
		&\ \Gamma_t\sumt\frac{3\alpha_i\gamma_i^2}{\Gamma_i}(\|\Delta_H^{2i}\|_*^2 + \|\Delta_H^{2i-1}\|_*^2) = \Gamma_t\sumt\frac{3\alpha_t\gamma_i^3}{\Gamma_t\gamma_t}(\|\Delta_H^{2i}\|_*^2 + \|\Delta_H^{2i-1}\|_*^2)
		=\frac{3\alpha_t}{\gamma_t}\sumt{\gamma_i^3}(\|\Delta_H^{2i}\|_*^2 + \|\Delta_H^{2i-1}\|_*^2),
	\end{align*}
	we conclude from \eqref{eqnEUC}, \eqref{tmp2} and Assumption \ref{itmVB} that
	\begin{align*}
		\E[\varepsilon\tn] 
		\leq&\ \Gamma_t\E[U_t] + \frac{3\alpha_t}{\gamma_t}\left[D^2+2\theta D^2 +  2\theta\Gamma_t\sumt{\gamma_i}\E[U_i] + \frac{6\alpha_t\sigma_H^2}{\gamma_t}\sumt{\gamma_i^3}\right] 
		\\
		\le &\ \frac{\alpha_t}{\gamma_t}C_t^2 +\frac{3\alpha_t}{\gamma_t}\left[(1+2\theta)D^2 + 2\theta\Gamma_t\sumt\frac{\alpha_i}{\Gamma_i} C_i^2 + \frac{6\alpha_t\sigma_H^2}{\gamma_t}\sumt{\gamma_i^3}\right].
	\end{align*}
	Finally, observing from \eqref{eqnCtheta} and \eqref{eqnGammaSpan} that $\ds\Gamma_t\sumt\frac{\alpha_i}{\Gamma_i}C_i^2\le C_t^2\Gamma_t\sumt\frac{\alpha_i}{\Gamma_i} = C_t^2$, we conclude \eqref{eqnEeps} from the above inequality.
	
\endproof
}

\vgap 
\section{Conclusion}
\label{secConclusion}
We present in this paper a novel accelerated mirror-prox (AMP) method for solving a
class of deterministic and stochastic variational inequality (VI) problems. The basic idea of
this algorithm is to incorporate a multi-step acceleration scheme into the mirror-prox method in \cite{nemirovski2005prox, juditsky2011solving}. For both the deterministic and stochastic VI, the AMP achieves the
optimal iteration complexity, not only in terms of its dependence on the number of  the iterations, but also on a variety of problem parameters. 
Moreover, the iteration cost of the AMP is comparable to, or even less than that of the mirror-prox method in that it saves one compuation of $\nabla G(\cdot)$.
To the best of our knowledge, this is the first algorithm with the optimal iteration complexity bounds for
solving the deterministic and stochastic VIs of type \eqref{eqnF}. Furthermore, we show that the developed AMP scheme
can deal with the situation when the feasible region is unbounded, as long as a strong solution of the VI exists.
In the unbounded case, we adopt the modified termination criterion employed by Monteiro and
Svaiter in solving monotone inclusion problem, and demonstrate that the rate of convergence of
AMP depends on the distance from the initial point to the set of strong solutions. Specially, in the unbounded case of the deterministic VI, the AMP scheme achieves the 
iteration complexity without requiring any knowledge on the distance from the initial point to the set of strong solutions.

\bibliography{all}{}
\bibliographystyle{siam_first_initial}

\end{document}